%% file: volumes.tex
\documentclass{amsart}
\pdfoutput=1

\usepackage{amsthm,amsmath,amssymb,latexsym,amsxtra,amscd,amssymb,color,mathtools,xspace}
\usepackage{mathtools}
\usepackage{units}
\usepackage{enumitem}
\usepackage{array}
\usepackage{booktabs}
\usepackage{calc}
\usepackage{pgf}
\usepackage{scalefnt}
\usepackage{listings}
\usepackage{multicol}
\usepackage{longtable}
\allowdisplaybreaks

\usepackage{bm}

\usepackage{placeins}
\usepackage[obeyspaces,hyphens,spaces]{url}
\usepackage{letltxmacro}

\usepackage{graphicx,color}
\graphicspath{{img/}}

\usepackage{tikz}
\usetikzlibrary{matrix,arrows,calc}
\usetikzlibrary[patterns]
\usetikzlibrary{decorations.markings}
\usetikzlibrary{decorations.pathmorphing,patterns}

\usepackage{csquotes}

\usepackage{enumitem}

\usepackage{amsthm}
\usepackage{hyperref}

\usepackage{aliascnt}

\input{macros}

\begin{document}

\title{Masur-Veech volume of the gothic locus}

\author{David Torres-Teigell}
\address{
Departamento de Matem\'{a}ticas, 
Universidad Aut\'{o}noma de Madrid,
28049 Madrid, Spain
}
\email{david.torres@uam.es}
\begin{abstract}
We calculate the Masur-Veech volume of the gothic locus $\G$ in the stratum $\cH(2^{3})$ of genus four. Our method is based on the use of the formulae for the Euler characteristics of gothic Teichm\"{u}ller curves to determine the number of lattice points of given area. We also use this method to recalculate the Masur-Veech volumes of the Prym loci $\cP_{3}\subset\cH(4)$ and $\cP_{4}\subset\cH(6)$ in genus three and four.
\end{abstract}

\maketitle
\tableofcontents


\input{sec_intro}

\input{sec_hodge}

\input{sec_hilbert}

\input{sec_teich}

\input{sec_number}

\input{sec_volumes}


%

\bibliographystyle{alpha}
\input{sources.bbl}


\end{document}

%% file: macros.tex
\def\equationautorefname~#1\null{(#1)\null}

\newtheorem{theorem}{Theorem}[section]

\newaliascnt{lemma}{theorem}
\newtheorem{lemma}[lemma]{Lemma}
\aliascntresetthe{lemma}

\newaliascnt{prop}{theorem}
\newtheorem{prop}[prop]{Proposition}
\aliascntresetthe{prop}

\newaliascnt{cor}{theorem}
\newtheorem{cor}[cor]{Corollary}
\aliascntresetthe{cor}

\newtheorem*{remark}{Remark}

\makeatletter%
\protected\def\W{\@ifnextchar[{\W@arg}{W_D}}
\def\W@arg[#1]{%
  \def\@D{D}%
  \def\c@W##1,##2,##3\@nil{\def\c@D{##1}\def\n@D{##2}}\c@W#1,,\@nil%
  \def\@W##1=##2=##3\@nil{\def\@DD{##1}\def\@DDD{##2}}\expandafter\@W\c@D==\@nil%
  \ifx\@D\@DD W_{\@DDD}%
    \expandafter\ifx\n@D\relax\relax\else(\n@D)\fi%
  \else W_D(#1)\fi%
}
\protected\def\parexp#1{\@ifnextchar^{(#1)}{#1}}
\protected\def\parind#1{\@ifnextchar_{(#1)}{#1}}

\def\defbb#1{\expandafter\def\csname b#1\endcsname{\mathbb{#1}}}
\def\defbm#1{\expandafter\def\csname bm#1\endcsname{\bm{#1}}}
\def\defcal#1{\expandafter\def\csname c#1\endcsname{\mathcal{#1}}}
\def\deffrak#1{\expandafter\def\csname frak#1\endcsname{\mathfrak{#1}}}
\def\defop#1{\expandafter\def\csname#1\endcsname{\operatorname{#1}}}
\def\defbold#1{\expandafter\def\csname bb#1\endcsname{\mathbf{#1}}}

\def\defcals#1{\@defcals#1\@nil}
\def\@defcals#1{\ifx#1\@nil\else\defcal{#1}\expandafter\@defcals\fi}
\def\deffraks#1{\@deffraks#1\@nil}
\def\@deffraks#1{\ifx#1\@nil\else\deffrak{#1}\expandafter\@deffraks\fi}
\def\defbbs#1{\@defbbs#1\@nil}
\def\@defbbs#1{\ifx#1\@nil\else\defbb{#1}\expandafter\@defbbs\fi}
\def\defbms#1{\@defbms#1\@nil}
\def\@defbms#1{\ifx#1\@nil\else\defbm{#1}\expandafter\@defbms\fi}
\def\defops#1{\@defops#1,\@nil}
\def\@defops#1,#2\@nil{\if\relax#1\relax\else\defop{#1}\fi\if\relax#2\relax\else\expandafter\@defops#2\@nil\fi}
\def\defbolds#1{\@defbolds#1\@nil}
\def\@defbolds#1{\ifx#1\@nil\else\defbold{#1}\expandafter\@defbolds\fi}
\makeatother%

\defbbs{ZHQCNPARDV}
\defbms{abvwz}
\defcals{OCFGHRMXPYZEVULAQNS}
\deffraks{abcdefopqgmB}
\defops{SL,mod,Spec,Re,Gal,Tr,tr,End,GL,Hom,Kern,Bild,PGL,Aut,id,Pic,div,Supp,supp,vol,area,Fix,Jac,im,Im,
PT,Pic,Id,Prym,diag,mult,tr,Sp,pr,Red,lcm,Per,RPer,Stab,rank}
\defbolds{uzvxeabcd}

\def\P{\cP}
\def\G{\cG}
\def\SS{S}
\def\a{\mathbf{a}}

\def\={\;=\;}  \def\+{\,+\,} 

\def\be{\begin{equation}}   \def\ee{\end{equation}}     
\def\bes{\begin{equation*}}    \def\ees{\end{equation*}}
\def\ba{\be\begin{aligned}} \def\ea{\end{aligned}\ee}   
\def\bas{\bes\begin{aligned}}  \def\eas{\end{aligned}\ees}

\newcommand{\R}[1][D]{R_{#1}}

\definecolor{lgray}{RGB}{214, 214, 214}
\definecolor{dgray}{RGB}{115, 115, 115}

%% file: sec_intro.tex
\section{Introduction}

The moduli space $\cH_{g}$ of flat surfaces of genus $g$ carries a remarkable $\SL(2,\bR)$-action preserving its natural stratification. By the famous theorem of Eskin, Mirzakhani and Mohammadi, orbit closures are very nice geometric objects called affine invariant manifolds, that is complex linear subspaces of a stratum locally defined by linear equations with real coefficients in period coordinates. The orbit closure of a generic flat surface is the whole stratum where it lives and, in fact, the list of proper affine invariant submanifolds of strata currently known is not very large.
\par
The main object of study in the present paper is the gothic locus $\G$, a four dimensional affine invariant manifold in the stratum $\cH(2^{3})$ of genus 4 discovered by McMullen, Mukamel and Wright in~\cite{MMW} (see~\autoref{subsec:prym_and_gothic} for definitions). This locus has the remarkable property that it contains one of the few known infinite families of geometrically primitive Teichm\"{u}ller curves. It is also closely related to the flex locus $F$, a totally geodesic, irreducible complex surface in $\cM_{1,3}$. In the present paper, we prove the following.

\begin{theorem}\label{thm:main} The Masur-Veech volume of the gothic locus $\G\subset \cH(2^{3})$ in genus 4 is
	\[\vol(\G) = \frac{13}{2^{7}\cdot 3^{5}} 
	\, \pi^{4}\,.\]
\end{theorem}

Using period coordinates, one can define the so-called Masur-Veech measure $\nu$ on a stratum $\cH(\bma)$ as the pullback of the Lebesgue measure on $H^{1}(X,Z(\omega);\bC)\cong\bC^{2g+n-1}$, the normalisation of which depends on the choice of a lattice in relative cohomology (see~\autoref{sec:hodge} and the references therein for discussions on normalisations). The measure of the whole stratum is obviously infinite, since one can always scale differentials. In order to solve this problem, one can restrict their attention to the hyperboloid $\cH_{1}(\bma)\subset \cH(\bma)$ formed by flat surfaces of $\area(X,\omega)=\tfrac{i}{2}\int_{X}\omega\wedge \overline{\omega}=1$. The measure of a set $B\subset \cH_{1}(\bma)$ is then defined to be the Masur-Veech measure of the cone $C(B)\subset \cH(\bma)$ over $B$. Masur and Veech proved independently~\cite{MasurMeasure,VeechMeasure} that this measure is finite. 
\par
More generally, if the coefficients of the equations defining an affine invariant manifold $\cN$ are rational numbers, one can define a measure on $\cN$ and $\cN_{1}$ in the same way as for the whole stratum (equivalently, one can consider the space of $\SL(2,\bR)$-invariant finite measures $\nu_{1}$ on $\cH_{1}(\bma)$ and consider their supports $\cN_{1}=\supp(\nu_{1})$, see~\cite{esmimo}).
\par
For strata $\cH(\bma)$ of flat surfaces almost everything is known: there are explicit formulae for the volumes, recursions to calculate them and large genus asymptotics (see~\cite{EO01, CMZ, Sauvaget, CMSZ}). The usual strategy for the computation of the volumes, due originally to Eskin and Masur and to Kontsevich and Zorich~\cite{ZorVol}, is based on the interpretation of the Lebesgue volume as the leading term of the number of lattice points in a cone as its radius tends to infinity. 
\par
In the case of strata $\cQ(\bmb)$ of finite-area quadratic differentials, there are precise formulae for the volumes in genus~0 and explicit values in small genera (see~\cite{EO06,AEZ,goujard}). Very recently, the authors of~\cite{vincent+} related the volume of the principal strata of quadratic differentials with the constant term of a family of polynomials and gave a recursive formula to compute them.

\

As far as we know, our result is the first example of computation of the Masur-Veech volume of a proper affine invariant submanifold of a stratum, apart from the loci of canonical double covers that can be naturally identified with strata of quadratic differentials. 
\par
Our method differs from the strategy of Eskin-Masur and Kontsevich-Zorich by a subtle but relevant twist: although we do calculate the precise number of minimal torus covers of given area and degree, we do not exhibit any of them.
\par
Instead, we use the formulae for the Euler characteristics of gothic Teichm\"{u}ller curves from~\cite{eulergothic} to give this number. The key observation is that the number of minimal torus covers $(X,\omega)\to (\bC/\bZ\oplus i\bZ,dz)$ of minimal area on an arithmetic Teichm\"{u}ller curve $C$ can be interpreted as the degree of the natural map $C\to M_{1,1}$ to the moduli space of elliptic curves, and this degree is given by the quotient of the Euler characteristics $\chi(C)/\chi(M_{1,1})$.
\par
The objective of exhibiting all minimal torus covers on $\cG$ is actually a challenging problem. In the case of full strata of abelian (resp. quadratic) differentials, one can take a combinatorial approach to the enumeration of square-tiled surfaces (resp. pillowcase covers) by determining all possible cylinder decompositions and counting the different square tilings (see for example~\cite{ZorVol, goujard}). However, one of the requirements to belong to the gothic locus is the existence of a non-Galois degree three map to an elliptic curve, and this condition does not have a direct interpretation in the flat geometry world. In order to enumerate the torus covers in $\G$ one would therefore need to translate this condition into a combinatorial feature that makes possible the determination of the possible cylinder decompositions. This approach, that would also allow the enumeration of cusps of gothic Teichm\"{u}ller curves, needs techniques beyond the scope of this paper.

\

The article is organised as follows. In~\autoref{sec:hodge} we review the basic aspects of the theory of flat surfaces, such as period coordinates on strata of abelian differentials, action of $\SL(2,\bR)$ and definition of the Masur-Veech measure both on strata and on affine invariant manifolds.~\autoref{sec:hilbert} is devoted to the study of Hilbert modular surfaces of square discriminant and abelian surfaces with real multiplication and the calculation of the induced polarisation on the abelian subvarieties generated by eigenforms. In~\autoref{sec:teich} we introduce the Prym and gothic loci and study the infinite families of Teichm\"{u}ller curves therein, giving formulae for their Euler characteristics and classifying the square-tiled surfaces they contain. In~\autoref{sec:numbertheory} we use Dirichlet series and modular forms to study the asymptotic behaviour of certain arithmetic functions related to the Euler characteristics of Prym and gothic Teichm\"{u}ller curves. Finally, in~\autoref{sec:volumes} we compute the volume of the gothic locus and the volumes of the minimal stratum in genus two and the two Prym loci in genus three and four. We also review how to link these last two values with the already known volumes of the corresponding strata of quadratic differentials.

\subsection*{Acknowledgements}
The author is grateful to Martin M\"{o}ller for proposing this project and for his continuous help during the development of the article. He also thanks Vincent Delecroix for very helpful discussions and the referee for a very careful reading and many suggestions. \\
The author was supported by the LOEWE-Schwerpunkt ``Uniformisierte Strukturen in Arithmetik und Geometrie'' and by the DFG-Projekt ``Classification of Teichm\"{u}ller curves MO 1884/2--1''.


%% file: sec_hodge.tex
\section{Masur-Veech volume and lattice points}\label{sec:hodge}

In this section we introduce the Masur-Veech volume on strata of abelian differentials and on affine invariant manifolds defined by rational equations in period coordinates. The normalisation of this measure depends on the choice of an appropriate lattice, the covolume of which is fixed to be one. The volume of an affine invariant manifold can then be calculated by counting the leading term of the number of lattice points in the ``hyperboloid'' of flat surfaces of area~$\le N$. This method, originally due to Kontsevich-Zorich and Eskin-Masur, and later completely developed by Eskin-Okounkov~\cite{EO01,EO06} using representation theory and quasimodular forms, has allowed the determination of general formulae for the volumes of different families of strata of abelian and quadratic differentials (see~\cite{AEZ}) and the calculation of these volumes in small genera (see~\cite{ZorVol} for the abelian case and~\cite{goujard} for the quadratic one). Here we follow this strategy to derive an explicit formula for the volume of affine invariant manifolds.

Before we state the main result of this section, let us introduce some of the objects needed. By a \emph{torus cover} $(X,\omega)\to (\bC/\Lambda,dz)$ we will mean a covering $\pi:X\to \bC/\Lambda$ such that $\pi^{*}dz=\omega$. The particular case $(X,\omega)\to(\bC/\Per(X,\omega),dz)$, whenever the lattice of periods $\Per(X,\omega)$ (see~\eqref{eq:Per_and_RPer}) is a sublattice of $\bZ\oplus i\bZ$, is called the \emph{minimal torus cover} of $(X,\omega)$. Note that, if we denote by $E$ the elliptic curve $E=\bC/\bZ\oplus i\bZ$, then for a torus cover $\pi:(X,\omega)\to (E,dz)$ one has $\area(X,\omega)=\deg(\pi)$ (see definitions below).

With the normalisations of~\autoref{subsec:MVmeasure}, the aforementioned method yields the following formula.

\begin{prop}\label{prop:volumeformula} 
Let $\cN$ be an affine invariant manifold defined by linear equations with rational coefficients in period coordinates. Assume that $\cN$ contains no rel deformation and normalise the Masur-Veech volume with respect to the lattice $\Lambda_{abs}^{\cN}$ of cohomology classes in $\cN$ with absolute periods in $\bZ\oplus i\bZ$. Define the sets
	\begin{align*}
	\begin{split}
	\cC_{d}(\cN) = \{(X,\omega)\in\cN \ :\ &\mbox{ $\pi:(X,\omega)\to (E,dz)$ torus cover of degree $d$} \}\,, \\ \label{eq:setsofsquaretiled}
	\cS_{d,m}(\cN) = \{(X,\omega)\in\cN \ :\ &
	\mbox{ $\pi:(X,\omega)\to (\bC/\Per(X,\omega),dz)$ minimal}
	\end{split}	\\
	&\mbox{ torus cover of degree $m$ and $\area(X,\omega)=d$} \}\,. \nonumber
	\end{align*}
 Then the volume of the hypersurface $\cN_{1}$ of flat surfaces of area~1 is given by
	\begin{align*}
	\vol(\cN_{1}) &= \lim_{D\to\infty} \frac{1}{D^{\dim_{\bC}\cN}}\sum_{d=1}^{D} \mid \cC_{d}(\cN) \mid \, = \\
	&= \lim_{D\to\infty} \frac{1}{D^{\dim_{\bC}\cN}}\sum_{d=1}^{D}\sum_{m|d} \sigma\bigl(\tfrac{d}{m}\bigr) \mid \cS_{m,m}(\cN) \mid \,,
	\end{align*}
where $\sigma(n)\coloneqq \sigma_{1}(n)$ denotes the sum of divisors of $n$.
\end{prop}

The assumption on the non-existence of rel deformations ensures that $\Lambda_{abs}^{\cN}$ is actually a lattice (see sections below). In~\autoref{sec:volumes} we will relate this counting problem with sums of Euler characteristics of arithmetic Teichm\"{u}ller curves in $\cN$.

\subsection{Strata of abelian differentials and period coordinates}

The \emph{Hodge bundle} $\cH_{g}$ over the moduli space $\cM_{g}$ is the moduli space of \emph{flat surfaces} of genus $g$, that is pairs $(X,\omega)$ where $X$ is a compact Riemann surface of genus $g$ and $\omega\in H^{0}(X,\Omega_{X})$ is a non-zero holomorphic 1-form, or abelian differential. A \emph{covering of flat surfaces} $\pi:(Y,\eta)\to(X,\omega)$ is simply a map between Riemann surfaces such that $\eta=\pi^{*}\omega$.

Given a flat surface $(X,\omega)$, the local charts given by $z \mapsto \int_{z_{0}}^{z}\omega$ around regular points $z_{0}$ and $z \mapsto \int_{z_{0}}^{z}\omega^{1/k}$ around zeroes $z_{0}$ of $\omega$ of order $k$ induce a translation structure on $X$, that is an atlas such that the transition functions on $X\setminus Z(\omega)$ are translations. This allows us to represent $(X,\omega)$ as a union of polygons in the plane whose sides are identified by translations, and it determines a flat metric on $X\setminus Z(\omega)$ given by the pullback of the euclidean metric on the plane. The flat area of $X$ is then given by the euclidean area of the polygon, which agrees with $\area(X,\omega)= \frac{i}{2}\int_{X} \omega\wedge\overline{\omega}$. Moreover, zeroes of order $k$ of the differential correspond to vertices of the polygon whose total angle sums up to $2\pi(k+1)$.

Now, the zero set $Z(\omega)$ of an abelian differential $\omega$ is a divisor of degree $2g-2$ on $X$ and one can therefore stratify $\cH_{g}$ according to the multiplicities of the zeroes of the differential, so that if $\bm{a}=(a_{1},\ldots,a_{n})$ is a partition of $2g-2$, the stratum $\cH(\bm{a})$ consists of flat surfaces $(X,\omega)$ such that $Z(\omega)=\sum a_{i}P_{i}$, for pairwise different points $P_{i}\in X$. The dimension of the stratum is $\dim_{\bC} \cH(\bba)=2g+n-1$. We will use exponential notation on the indices, for example $\cH(1^{2})\coloneqq \cH(1,1)$. 

Using the lattices $H^{1}(X;\bZ\oplus i\bZ)$ one can canonically identify nearby fibres of the Hodge bundle. This identification, equivalent to the Gauss-Manin connection, allows us to define the \emph{period coordinates} parametrising a stratum $\cH(\bm{a})$. In a small neighbourhood $U$ of a point $(X,\omega)\in\cH(\bma)$ one can define the \emph{period map}
	\[\begin{array}{ccl}
	U & \hookrightarrow & H^{1}(X,Z(\omega);\bC)\cong \bC^{2g+n-1} \\
	(Y,\eta) & \mapsto & \eta(\cdot)\,,
	\end{array}\]
where we identify $H^{1}(X,Z(\omega);\bC)=Hom(H_{1}(X,Z(\omega);\bZ),\bC)$ and the linear map $\eta(\cdot)$ is defined by $\eta(\gamma)=\int_\gamma \eta$ for $\gamma \in H_{1}(X,Z(\omega);\bZ)$. We also define the lattices of \emph{absolute and relative periods} of $(X,\omega)$ as
	\begin{align}\begin{split}
	\Per(X,\omega) &= \{ \omega(\gamma) \,:\, \gamma \in H_{1}(X;\bZ) \}\,, \\\label{eq:Per_and_RPer}
	\RPer(X,\omega) &= \{ \omega(\gamma) \,:\, \gamma \in H_{1}(X,Z(\omega);\bZ) \} \,.
	\end{split}\end{align}

Consider now the natural map $p:H^{1}(X,Z(\omega);\bC)\to H^{1}(X;\bC)$. A \emph{rel deformation} of $(X_{0},\omega_{0})\in\cH(\bba)$ is a path $\alpha(t)=(X_{t},\omega_{t})\in \cH(\bba)$, for $t\in[0,1]$, along which the lattice of absolute periods $\Per(X_{t},\omega_{t})$ stays the same, that is such that $p\circ\alpha$ is constant.

It is sometimes convenient to restrict ourselves to the hypersurface $\cH_{1}(\bma)\subset\cH(\bma)$ of flat surfaces of flat area~1. Note that, if $\{A_{j},B_{j}\,:\, i=1,\ldots,g\}$ are absolute periods with respect to a symplectic basis of homology, Riemann bilinear relations imply that $\area(X,\omega)= 
\frac{i}{2} \sum_{j} A_{j}\overline{B}_{j} - \overline{A}_{j}B_{j}$, and therefore $\cH_{1}(\bma)$ is locally a hyperboloid in period coordinates. The natural action of $\bC^{*}$ on $\cH(\bma)$ determines, for every flat surface $(X,\omega)$, a canonical point $(X,\lambda\omega)\in\cH_{1}(\bma)$ for some $\lambda\in\bR^{+}$. Note that $\area(X,\lambda\omega) = |\lambda|^{2}\cdot \area(X,\omega)$, for $\lambda\in\bC^{*}$.

\

In a similar way one can define strata of quadratic differentials. In this case the bundle $\cQ_{g}$ over the moduli space $\cM_{g}$ parametrises pairs $(Y,q)$ where $q$ is a non-zero quadratic differential on $Y$ with at most simple poles and not the square $q=\omega^{2}$ of an abelian differential $\omega \in H^{0}(Y,\Omega_{Y})$. One can consider local charts given by $z \mapsto \int_{z_{0}}^{z}\sqrt{q}$ as before in order to define a half-translation structure on $Y$, that is an atlas such that the transition functions on $Y\setminus Z(q)$ are given by $z \mapsto \pm z + c$. In this case the quadratic differential $(Y,q)$ can be represented as a union of polygons in the plane whose sides are identified by translations or translations composed with rotations of order 2, and the total angle around a vertex sums up to $\pi(k+2)$, where $k$ is the order of the zero (or pole) of $q$ at the corresponding point. Moreover, the condition on the orders of the poles ensures that the total area of $\int_{Y} |q|$ is finite.

Given a partition $\bmb=(b_{1},\ldots,b_{n})$ of $4g-4$, where the $b_{i}$ are non-zero integers $b_{i}\ge -1$, we denote by $\cQ(\bmb)$ the stratum of quadratic differentials such that $\div(q) = \sum b_{i}P_{i}$, which has dimension $\dim_{\bC} \cQ(\bbb)=2g+n-2$. In order to parametrise this stratum, consider the canonical double cover $\pi:X \to Y$ defined by the condition that $\pi^{*}q=\omega^{2}$ for an abelian differential $\omega\in H^{0}(X,\Omega_{X})$. If we denote by $\sigma: X\to X$ the involution that induces the map $\pi$, the abelian differential $\omega$ lies in the $-1$-eigenspace $H^{0}_{-}(X,\Omega_{X})$ of $\sigma^{*}$. Moreover, the multiplicities $\bma$ of $\omega$ are determined by those of $q$ in the following way: each zero of $q$ of even multiplicity $b_{i}$ yields two distinct zeroes of $\omega$ of multiplicities $b_{i}/2$, and each zero or pole of odd multiplicity $b_{i}$ yields a zero of $\omega$ of multiplicity $b_{i}+1$. 
The map $\cQ(\bmb) \to \cH(\bma)$ thus defined turns out to be an immersion, and the stratum $\cQ(\bbb)$ is then locally parametrised by the $-1$-eigenspace $H^{1}_{-}(X,Z(\omega);\bC)$ of $\sigma^{*}$ via the period map 
	\[\begin{array}{ccl}
	U & \hookrightarrow & H^{1}_{-}(X,Z(\omega);\bC)\cong \bC^{2g+n-2} \\
	(Y',q') & \mapsto & \omega'(\cdot)\,,
	\end{array}\]
where $\pi^{*}q'=(\omega')^{2}$.

Similarly to the abelian case, we will define the hypersurface $\cQ_{1}(\bmb)\subset\cQ(\bmb)$ of quadratic differentials of area~$1/2$. The reason for this normalisation comes from the fact that the canonical double cover $(X,\omega)$ of a quadratic differential $(Y,q)$ of area~$1/2$ has $\area(X,\omega)=1$.

\subsection{\texorpdfstring{$\SL(2,\bR)$}{SL(2,R)}-action and affine invariant manifolds}

There is a natural action of $\SL(2,\bR)$ and of $\GL^{+}(2,\bR)$ on the strata of abelian and quadratic differentials. In both cases, the action is induced by the affine action on the polygons in the plane representing the differential or, equivalently, by postcomposition of the charts of the translation or half-translation structure. Note that, in the abelian case, one can also construct $(X_{A},\omega_{A})\coloneqq A\cdot (X,\omega)$ where $A\in \GL^{+}(2,\bR)$ by considering the smooth differential
	\[\omega_{A} = \bigl( \begin{matrix} 1 & i \end{matrix} \bigr) \cdot A \cdot \Bigl( \begin{matrix}\, \Re \omega \, \\ \, \Im \omega \, \end{matrix} \Bigr)\]
on the differential manifold underlying $X$ and taking the unique Riemann surface structure $X_{A}$ with respect to which $\omega_{A}$ is holomorphic.

An \emph{affine invariant manifold} is a closed subset of a stratum $\cH(\bma)$ locally given by real linear equations in period coordinates. Since periods $\omega(\gamma)$ are precisely the holonomy vectors in $\bC\cong \bR^{2}$, that is the vectors representing the path $\gamma$ in the flat picture of $(X,\omega)$, the action of $\SL(2,\bR)$ in period coordinates is simply given by the affine action on $\bC^{N}\cong (\bR^{2})^{N}$. In particular, affine invariant manifolds are $\SL(2,\bR)$-invariant. The converse is also true: any $\GL^{+}(2,\bR)$-invariant closed analytic subspace is locally cut out by real linear equations (\cite[Prop. 1.2]{moellerlinear}). Furthermore, Eskin, Mirzakhani and Mohammadi proved the striking result that the closure of any $\GL^{+}(2,\bR)$-orbit is also an affine invariant manifold (\cite[Thm. 2.1]{esmimo}).

In the rare case that an orbit $\GL^{+}(2,\bR)(X,\omega)$ is already closed, it projects to an (algebraic) curve $C$ inside $\cM_g$, which is then called the \emph{Teichm\"{u}ller curve} generated by $(X,\omega)$. We will sometimes abuse the terminology and refer to the orbit itself as a Teichm\"{u}ller curve. Note that if $\pi:(Y,\pi^{*}\omega)\to (X,\omega)$ is a covering of flat surfaces, then the orbit $\GL^{+}(2,\bR)(Y,\pi^{*}\omega)$ is also closed. We say that a Teichm\"{u}ller curve is \emph{geometrically primitive} if it does not arise from such a covering construction. Not many families of geometrically primitive Teichm\"{u}ller curves are known, see e.g.~\cite{MMW} for a brief overview. Among them, McMullen and Calta independently constructed the \emph{Weierstra{\ss} curves} in genus $2$ (\cite{mcmtmchms,calta}) and McMullen later generalised this construction to the \emph{Prym-Teichm\"{u}ller curves} in genus $3$ and $4$ (\cite{mcmprym}). More recently, Eskin, McMullen, Mukamel and Wright found six exceptional orbit closures, two of which contain an infinite collection of Teichm\"{u}ller curves (\cite{MMW,EMMW}). We will be particularly interested in one of these families, the \emph{gothic Teichm\"{u}ller curves} in the stratum $\cH(2^{3})$. The explicit construction of these families of Teichm\"{u}ller curves will be carried out in~\autoref{sec:teich}.

Identifying the tangent space $T_{(X,\omega)}\cN$ with a subspace of $H^{1}(X,Z(\omega);\bC)$, the \emph{rank} of the affine invariant manifold $\cN$ is defined as the number $\tfrac{1}{2}\dim p(T_{(X,\omega)})$. Note that, whenever $\dim\cN = 2\rank(\cN)$, the kernel of $p$ restricted to $T_{(X,\omega)}\cN$ has dimension 0 and there are no non-trivial rel deformations inside $\cN$.

\subsection{Masur-Veech volume}\label{subsec:MVmeasure}

We consider the following loci in $H^{1}(X,Z(\omega);\bC)$
	\begin{align*}
	\Lambda_{rel}^{\cH(\bma)} &\coloneqq \{(Y,\eta)\in H^{1}(X,Z(\omega);\bC) \,:\, \RPer(Y,\eta) \subset \bZ\oplus i\bZ \}\,, \\
	\Lambda_{abs}^{\cH(\bma)} &\coloneqq \{(Y,\eta)\in H^{1}(X,Z(\omega);\bC) \,:\, \Per(Y,\eta) \subset \bZ\oplus i\bZ \}\,.
	\end{align*}
Note that $\Lambda_{rel}^{\cH(\bma)}$ agrees with $H^{1}(X,Z(\omega); \bZ\oplus i\bZ)$. On the other hand the locus $\Lambda_{abs}^{\cH(\bma)}$ does not define in general a lattice since rel deformations preserve it.

The most usual normalisation of the Masur-Veech volume uses the lattice $\Lambda_{rel}^{\cH(\bma)}$ (see~\cite[Appendix A]{DGZZ} or~\cite{goujard} for other discussions on the different normalisations). We also note here that points of these lattices represent certain torus covers, as we will see in the proof of~\autoref{prop:volumeformula}.

Take now the Lebesgue measure on $H^{1}(X,Z(\omega);\bC)\cong \bC^{2g+n-1}$, normalised so that the lattice $\Lambda_{rel}^{\cH(\bma)}$ has unit covolume. Its pullback via the period map is a well-defined and $\SL(2,\bR)$-invariant measure $\nu$ on $\cH(\bma)$ called the \emph{Masur-Veech measure}. This measure is obviously infinite, since one can always rescale $(X,\omega)$ inside $\cH(\bma)$. However, it allows us to define the following measure on the unit hyperboloid $\cH_{1}(\bma)$. For a subset $B\subset \cH_{1}(\bma)$ let us define the cone 
	\[C(B)=\{ (X,\lambda\omega)\in\cH(\bma) \,:\, (X,\omega)\in B\,,\ \lambda\in (0,1] \}\,.\]
We then define the volume of $B$ as $\vol(B) \coloneqq \nu(C(B))$ (other authors consider instead the measure defined by disintegration against the area -- see~\cite[\S 4.1]{AEZ} or~\eqref{eq:normalisation2_measure} below). Masur~\cite{MasurMeasure} and Veech~\cite{VeechMeasure} proved that this measure is actually finite, that is $\vol(\cH_{1}(\bma))<\infty$. 

One can similarly define a finite measure on the unit hyperboloid $\cQ_{1}(\bmb)$ by fixing a lattice in the space $H^{1}_{-}(X,Z(\omega);\bC)$ parametrising the stratum $\cQ(\bmb)$, see~\autoref{subsec:normalisations} for a quick review of the normalisations used in~\cite{AEZ,goujard}. 

\

Finally, let $\cN \subset \cH(\bma)$ be an affine invariant manifold defined by linear equations with rational coefficients in period coordinates and define the subspaces
	\ba\label{eq:latticenormalisation}
	\Lambda_{rel}^{\cN} &\coloneqq \Lambda_{rel}^{\cH(\bma)}\cap \cN = \{(Y,\eta)\in \cN\,:\, \RPer(Y,\eta) \subset \bZ\oplus i\bZ \}\,, \\
	\Lambda_{abs}^{\cN} &\coloneqq \Lambda_{abs}^{\cH(\bma)}\cap \cN = \{(Y,\eta)\in \cN\,:\, \Per(Y,\eta) \subset \bZ\oplus i\bZ \}\,,
	\ea
where we identify $\cN$ with its image under the period map. The rationality condition implies that $\Lambda_{rel}^{\cN}$ is a lattice in $\cN$.

If moreover $\cN$ does not admit rel deformations, that is $\dim\cN=2\rank(\cN)$, then $\Lambda_{abs}^{\cN}$ defines a lattice in $\cN$ as well. Since the gothic locus and the rest of affine invariant manifolds with which we will work satisfy this last condition, we will use this lattice $\Lambda_{abs}^{\cN}$ to normalise the Masur-Veech measure. 

In fact, completely analogously to the previous cases we can define the Masur-Veech measure $\nu_{\cN}$ on $\cN$ and the corresponding finite measure on $\cN_{1}$ with respect to the chosen lattice. We will sometimes abuse notation and write $\vol(\cN)$ instead of $\vol(\cN_{1})$ whenever the context is clear.

\begin{proof}[Proof of~\autoref{prop:volumeformula}]
By the definition of the Masur-Veech volume on $\cN_{1}$ with the lattice normalisation given by $\Lambda_{abs}^{\cN}$, one can approximate the volume $\vol(\cN_{1})$ by counting points in the lattice $\frac{1}{D_{0}}\Lambda_{abs}^{\cN}$ of area~$\le 1$, that is
	\[\vol(\cN_{1}) = \lim_{D_{0}\to\infty} \frac{1}{D_{0}^{\dim_{\bR}\cN}}\cdot \mid \bigl\{\tfrac{1}{D_{0}}\Lambda_{abs} \cap C(\cN_{1})) \bigr\} \mid \,.\]

Rescaling by $D_{0}$ and recalling that the area of $(X,\omega)$ grows quadratically under homotheties, we find that the following sets have the same cardinality:
	\[\mid \bigl\{\tfrac{1}{D_{0}}\Lambda_{abs}^{\cN} \cap C(\cN_{1})) \bigr\} \mid\ =\ \mid \bigl\{(X,\omega)\in \Lambda_{abs}^{\cN}\,:\, \area(X,\omega)\le D_{0}^{2} \bigr\}\mid \,.\]

Now, for each $(X,\omega)\in \Lambda_{abs}^{\cN}$, since $\Per(X,\omega)\subset \bZ\oplus i\bZ$ we have the following commutative diagram
	\begingroup

	\centering
	
	\begin{tikzpicture}
	  \matrix (m) [matrix of math nodes,row sep=2em,column sep=2em,minimum width=2em]
	  {
	     (X,\omega) & & (E,dz) \\
	     & (\bC/\Per(X,\omega),dz) & \\};
	  \path[-stealth]
	    (m-1-1) edge node [label={[label distance=-0.3cm]210:$\pi_{min}$}] {} (m-2-2.155)
	    (m-1-1) edge node [above] {$\pi$} (m-1-3)
	    (m-2-2.25) edge node [label={[label distance=-0.5cm]200:$j$}] {} (m-1-3)
	    ;
	\end{tikzpicture}

	\endgroup
\noindent given by $\pi(P)= \int_{P_{0}}^{P}\omega \bmod{\bZ\oplus i\bZ}$ and similarly for $\pi_{min}$. In particular, since $\area(E,dz)=1$, one has $\area(X,\omega) = \deg(\pi)$. 

Conversely, any degree $d$ covering $\pi:X\to E$ such that $\omega=\pi^{*}dz$ factorises as $\pi=j\circ \pi_{min}$ for the minimal torus covering $\pi_{min}:X\to\bC/\Per(X,\omega)$ and an isogeny $j: \bC/\Per(X,\omega) \to E$.

Writing $D=D_{0}^{2}$ proves the first equality of the formula.

\

In order to prove the second one, note that $\cC_{d}(\cN)=\bigcup_{m|d}\cS_{d,m}$. We claim that one has $\mid\cS_{d,m}(\cN)\mid\ = \ \sigma\bigl(\tfrac{d}{m}\bigr)\mid \cS_{m,m}(\cN)\mid$.

In fact, let us define the following function
	\[\begin{array}{cccl}
	\varphi : & \cS_{d,m}(\cN) & \to & \cS_{m,m}(\cN) \\
	 & (X,\omega) & \mapsto & M^{-1}\cdot(X, \omega)\,,
	\end{array}\]
where $M=M_{(X,\omega)}\in M(2,\bZ)\cap \GL^{+}(2,\bR)$ is the Hermite normal form associated to the lattice $\Per(X,\omega)$, that is the unique matrix
	\[M=\begin{pmatrix} a & s \\ 0 & c \end{pmatrix} \,,\quad 0<c\,,\ 0\le s< a\,,\ ac=[\bZ\oplus i\bZ : \Per(X,\omega)] = d/m\,,\]
such that $\Per(X,\omega)= a\bZ \oplus (s+ic)\bZ$ (see for example~\cite[\S 2.4]{Cohen}).

This map is well defined, since one has $\Per(M^{-1}\cdot(X,\omega))=M^{-1}\cdot \Per(X,\omega) = \bZ \oplus i\bZ$, and hence $\bC/\Per(M^{-1}\cdot(X,\omega))=E$. It is obviously surjective and, in fact, each $(X_{0},\omega_{0})\in \cS_{m,m}$ has as many preimages as sublattices of $\bZ\oplus i\bZ$ of index $d/m$ or, equivalently,
	\[\varphi^{-1}(\{(X_{0},\omega_{0})\}) = \left\{M\cdot(X_{0},\omega_{0})\,:\, M=\begin{pmatrix} a & s \\ 0 & \tfrac{d}{ma} \end{pmatrix}\, 0\le s < a \right\}\,,\]
which has cardinality $\sigma\bigl(\tfrac{d}{m}\bigr)$.
\end{proof}


%% file: sec_hilbert.tex
\section{Hilbert modular surfaces of square discriminant}\label{sec:hilbert}

Square-tiled surfaces generate (arithmetic) Teichm\"{u}ller curves for which part of the Jacobian admits real multiplication by a quadratic order $\cO_{d^{2}}$ of square discriminant. We will therefore be interested on Hilbert modular surfaces $X_{d^{2}}(\frakb_{r})$ (see the next two subsections for the definitions), which are the natural spaces parametrising such abelian varieties. In our case, the relevant part of the Jacobian will be an abelian surface endowed with a polarisation of type $(1,n)$. This will imply that the lattice defining it will be isomorphic as an $\cO_{d^{2}}$-module to the direct sum $\frakb_{r}\oplus \cO_{d^{2}}^{\vee}$, where $\frakb_{r}$ is an $\cO_{d^{2}}$-ideal of norm $n$. For each $d$, these ideals are determined by the choice of a positive integer $r|n$.

Abelian surfaces with real multiplication by $\cO_{d^{2}}$ always contain two natural elliptic curves generated by its eigenforms. These elliptic curves and their natural polarisation play an important role in the theory of arithmetic Teichm\"{u}ller curves. The following result determines their induced polarisation, which depends not just on the discriminant $d$, but also on the ideal $\frakb_{r}$ determining $X_{d^{2}}(\frakb_{r})$.

\begin{theorem}\label{prop:ellipticcurves}Let $[\bm\tau]\in X_{d^{2}}(\frakb_{r})$ and $T_{\bm\tau}$ be the associated $(1,n)$-polarised abelian surface with real multiplication by $\cO_{d^{2}}$ as in~\autoref{subsec:abeliansurfaces(1,n)}. Let $E_{1}, E_{2} \subset T_{\bm\tau}$ be the two elliptic curves generated by the two eigenforms for real multiplication $du_{1}$ and $du_{2}$. The restriction of the $(1,n)$-polarisation $\cL$ on $T_{\bm\tau}$ to $E_{1}$ and $E_{2}$ gives
	\[\cL|_{E_{1}} = \lcm(d,r)\cdot \cO_{E_{1}}(0)\quad\mbox{and}\quad 
	\cL|_{E_{2}} = \lcm(d,\tfrac{n}{r})\cdot \cO_{E_{2}}(0)\,. \]
\end{theorem}

This result will allow us to relate the discriminant of an arithmetic Teichm\"{u}ller curve and the degree of the minimal torus cover associated to any square-tiled surface on it.

\subsection{Hilbert modular surfaces.}\label{subsec:hilbert}

For any positive discriminant $D\equiv 0,1\mod 4$, write $D=b^2-4ac$ for some
$a,b,c \in\bZ$. The \emph{quadratic order of discriminant $D$}
is defined as $\cO_D=\bZ[T]/(aT^2+bT+c)$. 

In the case where $D=d^{2}$ for an integer $d>1$, one has $\cO_{d^{2}}=\bZ[T]/(T^{2}-dT)$ and this order is isomorphic to the subring
	\[ \cO_{d^{2}} = \{a=(a_{1},a_{2})\in\bZ\oplus\bZ \,:\, a_{1}\equiv a_{2} \bmod{d}\} \subset \bQ\oplus\bQ\]
with componentwise multiplication.

One can regard the pseudo-field $\bQ\oplus\bQ$ as an extension of $\bQ$ via the diagonal inclusion $\bQ\hookrightarrow \bQ\oplus\bQ$ and define a Galois conjugation $(a_{1},a_{2})^{\sigma}=(a_{2},a_{1})$. This allows us to define a norm $N(a)\coloneqq aa^{\sigma}=a_{1}a_{2}$ and a trace $\tr(a)\coloneqq a+a^{\sigma}=a_{1}+a_{2}$ as in the case of a field. The element $(d,-d)$ can therefore be interpreted as $\sqrt{D}$.

For any fractional ideal $\frakc\subset \bQ\oplus\bQ$, we denote by $\frakc^{\vee}$ the dual with respect to the trace pairing, i.e.\ $\frakc^{\vee}=\{x\in \bQ\oplus\bQ: \tr(x\frakc)\subset\bZ\}$. In particular, 
	\[\cO_{d^{2}}^{\vee}=\tfrac{1}{\sqrt{D}}\cO_{d^{2}}=\left\langle \tfrac{1}{d}(1,-1), (0,1)\right\rangle\,.\]
\par

Let $\frakb$ be an $\cO_{d^{2}}$-ideal. The $\cO_{d^{2}}$-module $\frakb\oplus\cO_{d^{2}}^{\vee}$ is preserved by the \emph{Hilbert modular group}
	\[\SL(\frakb\oplus\cO_{d^{2}}^{\vee}) \=
	\begin{pmatrix} \cO_{d^{2}} & \sqrt{D}\,\frakb \\ 
	\tfrac{1}{\sqrt{D}}\,\frakb^{-1} & \cO_{d^{2}} \end{pmatrix}\cap \SL_{2}(\bQ\oplus\bQ)\,,
	\]
which can be embedded into $\SL(2,\bQ)$ in two different ways induced by the two natural projections $\iota_{j}:\cO_{d^{2}}\to\bQ$. Associated with $\frakb$ we can construct the \emph{Hilbert modular surface}
	\[X_{d^{2}}(\frakb) \= \SL(\frakb\oplus\cO_{d^{2}}^{\vee})\backslash \bH^{2}\,,\]
where a matrix acts by M\"{o}bius transformations on each component via the two embeddings into $\SL(2,\bQ)$.

\subsection{Ideals of type $(1,n)$.}\label{subsec:ideals(1,n)}

The $\cO_{d^{2}}$-module $\frakb\oplus\cO_{d^{2}}^{\vee}$ can be equipped with a symplectic pairing of type $(1,1)$ given by
	\begin{equation}\label{eq:tracepairing}
	\langle (a,b), (x,y) \rangle  = \tr(ay-bx) = a_{1}y_{1}+a_{2}y_{2} - b_{1}x_{1} - b_{2}x_{2}\,.
	\end{equation}
In the coordinates $\cO_{d^{2}}=\langle (1,1),(0,d)\rangle$ and $\cO_{d^{2}}^{\vee}=\langle \tfrac{1}{d}(0,d),\tfrac{1}{d}(-1,1)\rangle$ this pairing is given by the matrix
	\[I=\left(\begin{matrix} 0 & J \\ -J^{T} & 0 \end{matrix}\right),\quad\mbox{ where } 
	J= \left(\begin{smallmatrix} 1 & 0 \\ d & 1 \end{smallmatrix}\right)\,.\]

Let $n$ be a square-free positive integer. We want to generalise this construction to find all rank-two $\cO_{d^{2}}$-modules with a symplectic pairing of signature $(1,n)$. By \cite{Bass62} such a lattice splits as a direct sum of $\cO_{d^{2}}$-modules, and is therefore isomorphic to $\frakb \oplus \cO_{d^{2}}^\vee$ for some $\cO_{d^{2}}$-ideal $\frakb$. This isomorphism can moreover be chosen so that the symplectic form is mapped to the trace pairing given by~\eqref{eq:tracepairing}. 
\par
In order to construct all such $\cO_{d^{2}}$-modules, let us define for each divisor $r$ of $n$ the ideal
	\[\frakb_{r} = \{(a_{1},a_{2})\in\bZ\oplus\bZ\,:\, 
	a_{1}\equiv a_{2}\bmod{d}\,,\ a_{1}\equiv 0\bmod{r}\,,\ a_{2}\equiv 0\bmod{n/r} \}\,.\]
The trace pairing on the associated $\cO_{d^{2}}$-module $\frakb_{r} \oplus \cO_{d^{2}}^\vee$ is of type $(r,n/r)$ or, equivalently, $(1,n)$. In this way one can generate all possible such $\cO_{d^{2}}$-modules.

\begin{lemma}\label{lem:ODmodules}Let $n$ be a square-free positive integer. Every $\cO_{d^{2}}$-module of rank 2 with a $(1,n)$-polarisation is isomorphic to $\frakb_{r} \oplus \cO_{d^{2}}^\vee$ for some divisor $r|n$. \\
Moreover, if $g=(d,n)$ and $\ell_{r}=\lcm(r,g)$ then $\frakb_{r}=\frakb_{s}$ if and only if $\ell_{r}=\ell_{s}$. In particular, $\frakb_{r}=\frakb_{\ell_{r}}$ for every $r|n$.
\end{lemma}

\begin{proof}First note that the ideal $\frakb$ must satisfy $\cO_{d^{2}}/\frakb\cong\bZ/n\bZ$. The corresponding map $\rho=\rho_{a}:\cO_{d^{2}}\to\bZ/n\bZ$ is determined by the image $\rho(0,d)=a$ that must satisfy $ad-a^{2}\equiv 0\bmod{n}$.

In the prime case $n=p$, the only possibilities are $a=0$ if $(d,p)=p$, and $a=0$ or $a=d$ if $(d,p)=1$. In both cases, $\ker \rho_{0} = \frakb_{p}$ and $\ker\rho_{d}= \frakb_{1}$. The general case follows since $n$ is square-free and the kernel of
	\[\begin{array}{ccc}
	\cO_{d^{2}} & \to & \prod \bZ/p_{i}\bZ \times \prod \bZ/p_{j}\bZ \\
	(0,d) & \mapsto & (0,\ldots,0,d,\ldots,d)
	\end{array}\]
is precisely $\frakb_{r}$ for $r=\prod p_{i}$.


The second claim follows from writing $g=g_{r}\cdot g_{n/r} \coloneqq \gcd(d,r)\cdot\gcd(d,n/r)$ and noting that both $g_{r}$ and $g_{n/r}$ necessarily divide each of the components of any element $(a_{1},a_{2})\in\frakb_{r}$.

\end{proof}

Note that Galois action on these ideals is given by $\frakb_{r}^{\sigma} = \frakb_{n/r}$. In particular one always has $\frakb_{r}^{\sigma}\neq \frakb_{r}$ unless $n|d$, in which case $\frakb_{1}=\frakb_{n}$ is the unique ideal of type $(1,n)$.

\begin{cor}\label{cor:ideals}Let $d>1$ and let $n$ be a square-free positive integer. The number of isomorphism classes of $\cO_{d^{2}}$-modules of rank 2 with a $(1,n)$-polarisation is $\sigma_{0}\bigl(\frac{n}{(d,n)}\bigr)$, where $\sigma_{0}$ counts the number of divisors. If $d\not\equiv 0 \bmod{n}$, each Galois orbit contains two elements.
\hfill$\square$
\end{cor}

We finish this section by explicitly calculating a basis for $\frakb_{r}$ in terms of $n$ and $d$. 

\begin{lemma}\label{lem:basisfrakb} Let $n$ be a square-free positive integer, $r|n$ and let $a,b\in\bZ$ be integers such that $ad+b\frac{n}{r}=(d,\frac{n}{r})$. Then
	\[\frakb_{r}=\left\langle (r \cdot(d,\tfrac{n}{r}), bn), (0,\lcm(d,\tfrac{n}{r})) \right\rangle_{\bZ}\,.\]
\end{lemma}

\begin{proof}The element $(0,\lcm(d,\tfrac{n}{r}))$ is primitive in $\frakb_{r}$ and the element $r \cdot (d,\tfrac{n}{r})\cdot(1,1)-ar\cdot(0,d)=(r\cdot (d,\tfrac{n}{r}),bn)$ belongs to $\frakb_{r}$. Since these elements are obtained from the usual generators of $\cO_{d^{2}}$ by the matrix 
	\[\left(\begin{matrix} r\cdot (d,\tfrac{n}{r}) & 0 \\ -ar & \frac{n/r}{(d,\tfrac{n}{r})} \end{matrix}\right)\]
of determinant $n$, the result follows.
\end{proof}

\subsection{Euler characteristics.}

The notion of Euler characteristic (of curves and of Hilbert modular surfaces)
refers throughout the text to orbifold Euler characteristics. The Euler characteristics of the classical Hilbert modular surfaces $X_D = X_D(\cO_D)$ were computed by Siegel in~\cite{siegel36}. A reference including also the square-discriminant case is \cite[Theorem~2.12]{bainbridge07}. The formula reads
\be\label{eq:chiXd2}
\chi(X_{d^{2}}) \= 2d^{3} \zeta_{\bQ\oplus\bQ}(-1) \sum_{r|d} \frac{\mu(r)}{r^2} = \frac{d^{3}}{72} \sum_{r|d} \frac{\mu(r)}{r^2}\,,
\ee
where $\mu$ is the M\"obius function and the zeta function $\zeta_{\bQ\oplus\bQ}$ is defined as the square of the usual Riemann zeta function. The following formula relates the Euler characteristics of the classical Hilbert modular surfaces $X_{d^{2}}$ and the ones in the case of a $(1,6)$-polarisation (\cite[Prop. 4.3]{eulergothic}):
	\be\label{eq:vol16}
\frac{\chi(X_{d^{2}}(\frakb_{r}))}{\chi(X_{d^{2}})} = \left\{\begin{array}{lcl}
1 \,,  & \text{if} & (6,d)=1 \\
3/2 \,, & \text{if} & (6,d)=2 \\
4/3 \,, & \text{if} & (6,d)=3 \\
2 \,,  & \text{if} & (6,d)=6. \\
\end{array} \right.
\ee

\subsection{Moduli of $(1,n)$-polarised abelian surfaces with real multiplication.}\label{subsec:abeliansurfaces(1,n)}

An abelian surface $T$ admits \emph{real multiplication} by $\cO_{D}$ if there exists an embedding $\cO_{D}\hookrightarrow \End(T)$ by self-adjoint endomorphisms. We will always assume that the action is \emph{proper}, in the sense that it cannot be extended to an action of a larger quadratic discriminant $\cO_{E}\supset \cO_{D}$.
\par
The different components of the moduli space of $(1,n)$-polarised abelian varieties with a choice of real multiplication by $\cO_D$ are parametrised by certain Hilbert modular surfaces (see~\cite[Chapter 7]{HvdG}). We will focus here in the square-discriminant case and follow similar lines as in~\cite[\S 4]{eulergothic} for the non-square case.
\par
In fact, if $(T = \bC^2/\Lambda,\cL)$ is an abelian variety with a $(1,n)$-polarisation $\cL$ and a choice of real multiplication by $\cO_{d^{2}}$, then $\Lambda$ is a rank-two $\cO_{d^{2}}$-module with symplectic pairing of signature $(1,n)$, hence it is isomorphic to some $\frakb_{r} \oplus \cO_{d^{2}}^\vee$ by~\autoref{lem:ODmodules}.
\par
Conversely, for any ideal $\frakb_{r}$ and
$\bm{\tau}=(\tau_{1},\tau_{2}) \in \bH^{2}$, we define the lattice
\be\label{eq:siegel}
\Lambda_{\frakb_{r},\bm{\tau}} \= \{(a_{1} + b_{1}\tau_1, a_{2}+b_{2} \tau_2)^T \,|\,\, 
a=(a_{1},a_{2}) \in \frakb_{r}\,,\ b=(b_{1},b_{2}) \in \cO_{d^{2}}^\vee \}.
\ee
The quotient $T_{\bm\tau}=\bC^{2}/\Lambda_{\frakb,\bm\tau}$ is an abelian surface 
with a $(1,n)$-polarisation (given by the trace pairing) and real multiplication
by $\cO_{d^{2}}$. The isomorphism class of $T_{\bm\tau}$ depends only on the image
of~$\bm\tau$ in $X_{D}(\frakb_{r})$. Note that the eigenforms $du_{1}$ and $du_{2}$ for real multiplication correspond to $(1,0)^{T}$ and $(0,1)^{T}$ in this basis.
\par
It follows from~\autoref{lem:ODmodules} and~\autoref{cor:ideals} that the locus of $(1,n)$-polarised abelian varieties with a choice of real multiplication~by $\cO_{d^{2}}$ has $\sigma_{0}(n/(d,n))$ components, each of these components being parametrised by a Hilbert modular surface $X_{D}(\frakb_{r})$. This allows us to characterise this locus in the case that we are studying. For the rest of the article we will keep the following explicit choices for the representatives $r|6$.

\begin{prop}\label{prop:Xcomponents}
The moduli space of $(1,6)$-polarised abelian surfaces with a choice of real
multiplication by $\cO_{d^{2}}$ consists of the following Hilbert modular surfaces
	\begin{itemize}
	\item $X_{d^{2}}(\frakb_{1})$ for $d \equiv 0\bmod{6}$;
	\item $X_{d^{2}}(\frakb_{1})\cup X_{d^{2}}(\frakb_{2})$ for $d \equiv 3\bmod{6}$,
	\item $X_{d^{2}}(\frakb_{1})\cup X_{d^{2}}(\frakb_{3})$ for $d \equiv 2,4\bmod{6}$, and
	\item $X_{d^{2}}(\frakb_{1})\cup X_{d^{2}}(\frakb_{2})\cup X_{d^{2}}(\frakb_{3})\cup X_{d^{2}}(\frakb_{6})$ for $d \equiv 1,5\bmod{6}$.
	\end{itemize}
\hfill$\square$
\end{prop}
\par
Note, however, that the locus of real multiplication in $\cA_{2,(1,6)}$
has in general fewer components than the moduli space of abelian surfaces
with a chosen real multiplication by~$\cO_{d^{2}}$. In fact, the
abelian varieties parametrised by $X_{d^{2}}(\frakb)$ and by $X_{d^{2}}(\frakb^\sigma)$
map to the same subsurface in $\cA_{2,(1,6)}$.

We can now prove the main result of the section. 

\begin{proof}[Proof of~\autoref{prop:ellipticcurves}] The covering map $\pi:X\to E$ associated to a square-tiled surface $(X,\omega)$ induces a map $\pi^{*}:E \to \Jac X$, which is injective if $\pi$ is primitive. Since in this situation $\omega=\pi^{*}dz$, we will be interested in the case where an eigenform for real multiplication by $\cO_{d^{2}}$ generates an abelian subvariety $E\subset \Prym X$.

Let $X_{D}(\frakb_{r})$ parametrise a component of the moduli space of $(1,n)$-polarised abelian varieties with a choice of real multiplication by $\cO_{d^{2}}$ as above. Recall the forgetful map
	\[\begin{array}{ccc} 
	X_{D}(\frakb_{r}) & \to & \cA_{2,(1,n)} \\
	\left[ \bm\tau\right] & \mapsto & T_{\bm\tau}
	\end{array}
	\]
to the moduli space of $(1,n)$-polarised abelian varieties as in~\eqref{eq:siegel}.

By~\autoref{lem:basisfrakb}, the period matrix of the abelian variety $T_{\bm\tau}$ in the eigenform basis is given by
	\[\Pi=
	\begin{pmatrix}
	r \cdot (d,\tfrac{n}{r}) & 0 & 0 & -\frac{1}{d}\tau_{1} \\
	bn & \lcm(d,\tfrac{n}{r}) & \tau_{2} & \frac{1}{d}\tau_{2}
	\end{pmatrix}
	\,,\]
for some integer $b$ coprime to $d$.


Now note that the elliptic curves $E_{1}$ and $E_{2}$ are generated by the eigenforms for real multiplication, that is
	\[E_{1} = \frac{\bC du_{1}}{\Lambda_{1}} \subset T_{\bm\tau}\quad\mbox{ and }\quad
	E_{2} = \frac{\bC du_{2}}{\Lambda_{2}} \subset T_{\bm\tau}\,,\]

where
	\begin{align*}
	\Lambda_{1} &= \langle du_{1}\rangle \cap \Lambda_{\frakb_{r},\bm\tau} = \left\langle \left(\begin{array}{c} \lcm(d,r) \\ 0 \end{array}\right) , \left(\begin{array}{c} \tau_{1} \\ 0 \end{array}\right) \right\rangle_{\bZ} \\
	\Lambda_{2} &= \langle du_{2}\rangle \cap \Lambda_{\frakb_{r},\bm\tau} = \left\langle \left(\begin{array}{c} 0 \\ \lcm(d,\frac{n}{r}) \end{array}\right) , \left(\begin{array}{c} 0 \\ \tau_{2} \end{array}\right) \right\rangle_{\bZ}
	\end{align*}

The result follows from restricting the symplectic pairing~\eqref{eq:tracepairing} to these sublattices.
\end{proof}

\subsection{Boundary components of Hilbert modular surfaces}

Recall the ideals
	\begin{align*}
	\frakb_{r} &\coloneqq \{(a_{1},a_{2})\in \bZ^{2}\,:\, a_{1}\equiv a_{2}\bmod{d}\,,\ a_{1}\equiv 0\bmod{r}\,, a_{2}\equiv 0\bmod{\tfrac{6}{r}}\} \,, \\
	\frakb_{6/r} &\coloneqq \{(a_{1},a_{2})\in \bZ^{2}\,:\, a_{1}\equiv a_{2}\bmod{d}\,,\ a_{1}\equiv 0\bmod{\tfrac{6}{r}}\,, a_{2}\equiv 0\bmod{r}\} = \frakb_{r}^{\sigma} 
	\end{align*}
and consider the Hilbert modular groups
	\begin{align*}
	\SL(\frakb_{r}\oplus\cO_{d^{2}}^{\vee}) &=
	\begin{pmatrix} \cO_{d^{2}} & \sqrt{D}\,\frakb_{r} \\ 
	\tfrac{1}{6\sqrt{D}}\,\frakb_{6/r} & \cO_{d^{2}} \end{pmatrix}\cap \SL_{2}(\bQ\oplus\bQ)\,, \\
	\SL(\frakb_{r}\oplus\cO_{d^{2}}) &=
	\begin{pmatrix} \cO_{d^{2}} & \frakb_{r} \\ 
	\tfrac{1}{6}\,\frakb_{6/r} & \cO_{d^{2}} \end{pmatrix}\cap \SL_{2}(\bQ\oplus\bQ)\,,
	\end{align*}
where $\sqrt{D}=(d,-d)$. These two groups are conjugate by $\left(\begin{smallmatrix} 1 & 0 \\ 0 & \sqrt{D} \end{smallmatrix}\right)$ and the corresponding Hilbert modular surfaces
	\[\SL(\frakb_{r}\oplus\cO_{d^{2}}^{\vee})\backslash \bH\times\bH\quad\mbox{ and }\quad 
	\SL(\frakb_{r}\oplus\cO_{d^{2}})\backslash \bH\times(-\bH)\]
are isomorphic under $(\tau_{1},\tau_{2})\mapsto (d\tau_{1},-d\tau_{2})$. We will therefore sometimes use the latter to simplify notation.

The Bairy-Borel compactification of the Hilbert modular surface $X_{D}(\frakb)$ is given by
	\[\widehat{X}_{D}(\frakb) = \SL(\frakb\oplus \cO_{D}^{\vee})\backslash (\bH\times\bH)_{D}\]
where
	\[(\bH\times\bH)_{D} = \left\{\begin{array}{ll}
	(\bH\times\bH) \cup \bP^{1}(\bQ(\sqrt{D}))\,, & \mbox{if $D\neq d^{2}$,}\\
	(\bH\cup\bP^{1}(\bQ))\times (\bH\cup\bP^{1}(\bQ)) \,, & \mbox{if $D=d^{2}$,}\\
	\end{array}\right.\]
and the action of $\SL(\frakb\oplus \cO_{D}^{\vee})$ extends in the natural way (see~\cite{bainbridge07,vdG} for background). In the case of quadratic discriminants $D=d^{2}$, the Baily-Borel compactification $\widehat{X}_{d^{2}}(\frakb_{r})$ has just orbifold singularities and its boundary is formed by the curves
	\[
	S_{d^{2}}^{1}(\frakb_{r}) = \bigcup_{p\in\bP^{1}(\bQ)} \overline{\bH\times \{p\}}\quad \mbox{ and } \quad
	S_{d^{2}}^{2}(\frakb_{r}) = \bigcup_{p\in\bP^{1}(\bQ)} \overline{\{p\}\times\bH}
	\]
in $\widehat{X}_{d^{2}}(\frakb_{r})$ and the set of cusps. The number of components of the curves $S_{d^{2}}^{j}(\frakb_{r})$ depends on $r$ and on the value of $d$, but we will see that it is bounded.

In order to describe the connected components of $S_{d^{2}}^{j}(\frakb_{r})$ we need to study the stabilisers of cusps in the first and second components of $\bH^{2}$. The different cases depend heavily on $d \bmod{6}$ and a precise description of the geometry of these connected components exceeds the objectives of this paper. The following estimate is enough for our purposes.

\begin{lemma}\label{lem:componentsboundary} The number of irreducible components of $S_{d^{2}}^{j}(\frakb_{r})$ is at most 12. Each of these irreducible components is isomorphic to a curve $K\backslash \bH$, where $K$ is a Fuchsian group containing $\Gamma_{1}(N)$ for $N=2160\, d$.
\end{lemma}

\begin{proof} The number of components of $S_{d^{2}}^{2}(\frakb_{r})$ is given by the number of orbits of $\bP^{1}(\bQ)\times \bH$ under $\SL(\frakb_{r}\oplus\cO_{d^{2}})$. By solving for its second coordinate, one can prove that the action of $\SL(\frakb_{r}\oplus\cO_{d^{2}})$ on the first coordinate includes the action of the group $A^{-1} \Gamma_{0}(6r)A$, where $A=\left(\begin{smallmatrix}0 & 1 \\ -1 & 0 \end{smallmatrix}\right)$. In particular the number of components of $S_{d^{2}}^{2}(\frakb_{r})$ is bounded by the number of cusps of $\Gamma_{0}(6r)$, which is at most $C(\Gamma_{0}(36))=12$.

We now want to prove that the stabiliser of each component $\overline{\{p\}\times\bH}$ contains a subgroup isomorphic to $\Gamma_{1}(2160\,d)$. This is trivial for the cusp at infinity by considering $B^{-1} \Gamma_{1}(36d)B$, with $B=\left(\begin{smallmatrix}1 & 0 \\ 0 & 6 \end{smallmatrix}\right)$. In general, a case by case calculation shows that, apart from 0 and $\infty$, all cusps of the groups $A^{-1} \Gamma_{0}(6r)A$ have representatives of the form $p=-\tfrac{36}{m}$ for some $m|360$. In particular, the matrix $M = \left(\begin{smallmatrix}1 & 0 \\ m/36 & 1 \end{smallmatrix}\right)$ sends this cusp to $\infty$ and, therefore, studying the cusp at $(p,p)$ of the Hilbert modular group $\SL(\frakb_{r}\oplus\cO_{d^{2}})$ is equivalent to studying the cusp at $(\infty,\infty)$ of the Hilbert modular group
	\[M\cdot \SL(\frakb_{r}\oplus\cO_{d^{2}})\cdot M^{-1} = \SL(\frakb_{r}\fraka^{-1}\oplus\cO_{d^{2}})\,,\quad \mbox{ where } \fraka=\tfrac{m}{36}\frakb_{r} + \cO_{d^{2}}\,. \]

Using the inclusions $\tfrac{m}{36}\frakb_{r} \subset \fraka \subset \tfrac{1}{36}\cO_{d^{2}}$, it is a direct computation to show that, for $C= \left(\begin{smallmatrix} 1 & 0 \\ 0 & 216 \end{smallmatrix}\right)$, the stabiliser $\Stab_{\{p\}\times\bH}\SL(\frakb_{r}\fraka^{-1}\oplus\cO_{d^{2}})$ contains a subgroup isomorphic to $C^{-1}\Gamma_{1}(N)C$. In fact, for each matrix $\left(\begin{smallmatrix} a & 216b \\ \frac{360cd}{36} & e \end{smallmatrix}\right) \in C^{-1}\Gamma_{1}(N)C$ the pair
	\[
	\left( 
	\begin{pmatrix} 1 & 216b \\ 0 & 1 \end{pmatrix}\,,
	\begin{pmatrix} a & 216b \\ \frac{360cd}{36} & e \end{pmatrix}
	\right) \in \begin{pmatrix} \cO_{d^{2}} & 36\frakb_{r} \\ \frac{m}{36}\cO_{d^{2}} & \cO_{d^{2}} \end{pmatrix} \cap \SL_{2}(\bQ\oplus \bQ)\,.
	\]
\end{proof}

%% file: sec_teich.tex
\section{Affine invariant manifolds and Teichm\"{u}ller curves}\label{sec:teich}

Because of the formula for the volume in~\autoref{prop:volumeformula}, we will be interested in counting the number $\mid \cS_{m,m}(\cN)\mid$ of minimal torus covers with fixed degree and area~$m$ in the different affine invariant manifolds. It is clear that they all belong to certain Teichm\"{u}ller curves of square discriminant $D=d^{2}$, but we need the precise relation between $m$ and $d$ in order to be able to apply the formulae for the Euler characteristics.

More precisely, in the minimal stratum $\cH(2)$ of genus 2, McMullen proved in~\cite[\S 6]{mcmTCingenustwo} that any minimal torus cover $\pi_{min}:(X,\omega)\to(\bC/\Per(X,\omega),dz)$ of degree $m$ and $\area(X,\omega)=m$ belongs to $W_{m^{2}}(2)$. In~\cite[Proposition B.1]{lanneaunguyen} and~\cite[Proposition 4.2]{LNPrym4} Lanneau and Nguyen proved the equivalent result for Prym-Teichm\"{u}ller curves in genus 3 and 4 using cylinder decompositions. 

\

Here we use \autoref{prop:ellipticcurves} to classify minimal torus covers in the gothic locus. This method extends easily to any affine invariant manifold and depends only on the polarisation on the part of the Jacobian admitting real multiplication. We follow the definitions and notation introduced in Sections~\ref{subsec:prym_and_gothic} and~\ref{subsec:euler} below.


\begin{theorem}\label{thm:areadiscgothic} For each $m>0$, the set of minimal torus covers of degree and area $m$ in the gothic locus $\cG$ distributes among different Teichm\"{u}ller curves in the following way:
	\[
	\cS_{m,m}(\G) = \left\{
	\begin{array}{l}
	\cS_{m,m}(G^{1}_{m^{2}}) \cup \cS_{m,m}(G^{2}_{(m/2)^{2}}) \cup \cS_{m,m}(G^{3}_{(m/3)^{2}}) \cup \cS_{m,m}(G^{6}_{(m/6)^{2}}) \,, \hspace{0.9cm} \\
	\omit$\displaystyle\hfill \mbox{if $m\equiv 6,30 \bmod{36}$,}$ \\[4pt]
	\cS_{m,m}(G^{1}_{m^{2}}) \cup \cS_{m,m}(G^{2}_{(m/2)^{2}}) \,, \\
	\omit$\displaystyle\hfill\text{if $m\not\equiv 6,30 \bmod{36}$ and $m\equiv 2 \bmod{4}$,}$ \\[4pt]
	\cS_{m,m}(G^{1}_{m^{2}}) \cup \cS_{m,m}(G^{3}_{(m/3)^{2}}) \,, \\
	\omit$\displaystyle\hfill \mbox{if $m\not\equiv 6,30 \bmod{36}$ and $m\equiv 3,6 \bmod{9}$,}$ \\[4pt]
	\cS_{m,m}(G^{1}_{m^{2}}) \,, \\
	\omit$\displaystyle\hfill \mbox{otherwise.}$
	\end{array}\right.
	\]

\end{theorem}

\

\subsection{The Prym and gothic loci}\label{subsec:prym_and_gothic}

Given an involution $J:X\to X$ of a compact Riemann surface $X$, we will denote by $\pi_{J}:X\to X/J$ the quotient map and by $H^{0}_{-}\coloneqq H^{0}_{-}(X,\Omega_{X})$ the $-1$-eigenspace of $J^{*}$. We will also say that a map $\pi:X\to B$ to an elliptic curve $B$ is an \emph{odd map} (with respect to $J$) if there exists an involution $j:B\to B$ such that $\pi_{B} \circ J= j\circ \pi_{B}$.

The \emph{Prym loci} in genus 3 and 4 and the \emph{gothic locus} in genus 4 are affine invariant manifolds of (complex) dimension four containing infinite families of geometrically primitive Teichm\"{u}ller curves (see~\cite{mcmprym, MMW}). They are defined in the following way
	\begin{align*}
	\cP_{3} &= \{(X,\omega)\in\cH(4)\,:\, &&\hspace{-0.4cm}\exists\ J:X\to X \mbox{ involution, } g(X/J)=1\,,\ \omega\in H^{0}_{-} \}\,, \\
	\cP_{4} &= \{(X,\omega)\in\cH(6)\,:\, &&\hspace{-0.4cm}\exists\ J:X\to X \mbox{ involution, } g(X/J)=2\,,\ \omega\in H^{0}_{-} \}\,, \\
	\cG &= \{(X,\omega)\in\cH(2^{3})\,:\, &&\hspace{-0.4cm}\exists\ J:X\to X \mbox{ involution, } g(X/J)=1\,,\ \omega\in H^{0}_{-}\,, \\
	&&&\hspace{-0.4cm}\exists\ \pi_{B}:X\to B \mbox{ odd map, }\deg(\pi)=3\,,\ |\pi(Z(\omega))|=1\}\,.
	\end{align*}

They are all defined by linear equations with integral coefficients in period coordinates (see~\cite[(9.2)]{MMW} for the gothic case).

For a surface $(X,\omega)$ in one of the Prym loci, the space of holomorphic differentials decomposes as $H^{0}(X,\Omega_{X})=\pi_{J}^{*}(H^{0}(X/J,\Omega_{X/J}))\oplus H^{0}_{-}(X,\Omega_{X})$. The corresponding lattices in homology $H_{1}^{-}(X;\bZ) \coloneqq H_{1}(X;\bZ) \cap (H^{0}_{-}(X,\Omega_{X}))^{*}$ carry a polarisation of type $(1,2)$ and $(2,2)$ for $\cP_{3}$ and $\cP_{4}$, respectively. One defines the \emph{Prym variety} of $(X,\omega)$ as the abelian subvariety $\Prym X\coloneqq (H^{0}_{-}(X,\Omega_{X}))^{*}/H_{1}^{-}(X;\bZ)$ of the Jacobian $\Jac X$.

As for  a surface $(X,\omega)$ in the gothic locus, the space of differentials decomposes as $H^{0}(X,\Omega_{X})=\pi_{J}^{*}(H^{0}(X/J,\Omega_{X/J}))\oplus \pi_{B}^{*}(H^{0}(B,\Omega_{B})) \oplus H^{0}_{\G}(X,\Omega_{X})$, for some subspace $H^{0}_{\G}(X,\Omega_{X})$ containing $\omega$. The lattice $H_{1}^{\G}(X;\bZ) \coloneqq H_{1}(X;\bZ) \cap (H^{0}_{\G}(X,\Omega_{X}))^{*}$ carries this time a polarisation of type $(1,6)$. We define equivalently the \emph{Prym variety} of $(X,\omega)$ as the abelian surface $\Prym X\coloneqq (H^{0}_{\G}(X,\Omega_{X}))^{*}/H_{1}^{\G}(X;\bZ)$ inside the Jacobian $\Jac X$.

We will also consider the stratum $\cH(2)$ of genus two forms $(X,\omega)$ with a single zero and abuse notation by writing $\Prym X=\Jac X$ in this case.

In each of these four affine invariant manifolds and for every discriminant $D>0$ the subspace 
	\begin{align*}
	\cN_{D}=\{(X,\omega)\in\cN\,:\ & \Prym X \mbox{ admits real multiplication by $\cO_{D}$} \\
	& \mbox{ with $\omega$ as an eigenform}\}
	\end{align*}
of eigenforms for real multiplication forms a closed $\GL^{+}(2,\bR)$-orbit. The corresponding Teichm\"{u}ller curves are denoted by $W_{D}(2)\subset\cH(2)$, $W_{D}(4)\subset \cP_{3}$, $W_{D}(6)\subset \cP_{4}$ and $G_{D}\subset \G$, respectively.

\subsection{Euler characteristics of Teichm\"{u}ller curves}\label{subsec:euler}

The Euler characteristics of the infinite families of Teichm\"{u}ller curves just defined have been computed by various authors. 

In genus 2, Bainbridge computed the Euler characteristics of non-arithmetic Weierstra{\ss}-Teichm\"{u}ller curves $W_{D}(2)$. Contrary to the Prym and gothic case, there is also an explicit formula for the Euler characteristics of arithmetic Weierstra{\ss}-Teichm\"{u}ller curves, given by Eskin, Masur and Schmoll. This allows us to give an easy asymptotic formula for the Euler characteristics of Hilbert modular surfaces with square discriminant.

\begin{theorem}[\cite{bainbridge07,EMS}]\label{thm:eulerg2} The Euler characteristics of the Weierstra{\ss} Teichm\"{u}ller curves in genus 2 are given by
	\begin{align*}
	\chi(W_{D}(2)) &= -\tfrac{9}{2} \chi(X_{D})\,,\quad\mbox{for $D>4$ a non-square discriminant,}\\
	\chi(W_{d^{2}}(2)) &= -\tfrac{d^{2}(d-2)}{16} \sum_{r|d}\frac{\mu(r)}{r^{2}}\,,\quad\mbox{for $d>1$.}
	\end{align*}
\end{theorem}

The Euler characteristics of the Prym-Teichm\"{u}ller families $W_{D}(4)$ and $W_{D}(6)$ in genus 3 and 4 were calculated by M\"{o}ller~\cite{moellerprym}, and the number of connected components was given by Lanneau and Nguyen~\cite{lanneaunguyen,LNPrym4}. The situations in genus 3 and in genus 4 differ drastically. Whereas Prym-Teichm\"{u}ller curves $W_{D}(6)$ in genus 4 behave basically like the Weierstra{\ss} curves in genus 2, the formulae for the Euler characteristic and the number of connected components of the ones in genus 3 depend on the arithmetic structure of $D$.

\begin{theorem}[\cite{moellerprym,lanneaunguyen,LNPrym4}]\label{thm:eulerprym} Let $D=f^{2}D_{0}>4$ be a non-square discriminant with conductor $f$.

\noindent $\bullet$ The Prym-Teichm\"{u}ller curve $W_{D}(6)$ in genus 4 is irreducible and its Euler characteristic is given by
	\[
	\chi(W_{D}(6)) = - 7 \chi(X_{D})\,,\quad\mbox{for $D>4$ a non-square discriminant.}
	\]
\noindent $\bullet$ The Prym-Teichm\"{u}ller curve $W_{D}(4)$ in genus 3 is empty if $D\equiv 5\bmod{8}$, has one irreducible component $W_{D}^{1}(4)$ if $D\equiv 0,4\bmod{8}$ and two irreducible components $W_{D}^{1}(4) \cup W_{D}^{2}(4)$ if $D\equiv 1\bmod{8}$. The Euler characteristic of each component $W_{D}^{j}(4)$ is given by
	\[
	\chi(W_{D}^{j}(4)) = 
	\left\{
	\begin{array}{ll}
	-\tfrac{5}{2}\, \chi(X_{D})\,, & \mbox{if $2\nmid f$,} \\
	-\tfrac{15}{4}\, \chi(X_{D})\,, & \mbox{if $2 \mid f$.}
	\end{array}
	\right.
	\]
\end{theorem}


Finally, the Euler characteristics of the Gothic-Teichm\"{u}ller curves were recently computed in a joint work with M\"{o}ller.

\begin{theorem}[\cite{eulergothic}]\label{thm:eulergothic} Let $D=f^{2}D_{0}>5$ be a non-square discriminant with conductor $f$. The Gothic-Teichm\"{u}ller curve $G_D$ is non-empty if and only if $D \equiv 0,1,4,9,12,16\bmod{24}$. \\ 
In this case, $G_{D}$ consists of several (perhaps still reducible) components. The number of such components agrees with $c_{D}$ the number of ideals of norm $6$ in $\cO_{D}$, that is one if $D \equiv 0,12\bmod{24}$, two components if $D \equiv 4,9,16\bmod{24}$ and four components if $D \equiv 1\bmod{24}$. For a fixed $D$ the Euler characteristics of all these components agree and are equal to
	\[\chi(G_D^{j}) =  
	\left\{\begin{array}{lcl}
	-\tfrac32  \, \chi(X_D) - 2\, \chi(\R^{j}) \,, & \text{if} & (6,f)=1 \\
	-\tfrac94  \, \chi(X_D) - 2\, \chi(\R^{j}) \,, & \text{if} & (6,f)=2 \\
	-2         \, \chi(X_D) - 2\, \chi(\R^{j}) \,, & \text{if} & (6,f)=3 \\
	-3         \, \chi(X_D) - 2\, \chi(\R^{j}) \,, & \text{if} & (6,f)=6. \\
	\end{array} \right.\]
\end{theorem}

Here $\R^{j}$ denotes the \emph{$(2,3)$-reducible locus}, which is the pullback to $X_{D}(\frakb_{j})$ of the locus inside the moduli space $\cA_{2,(2,3)}$ of $(2,3)$-polarized abelian surfaces consisting of products $E_{1}\times E_{2}$ of elliptic curves with the natural
$(2,3)$-polarization $2\,p_{1}^{*}\cO_{E_{1}}(0)\otimes 3\,p_{2}^{*}\cO_{E_{2}}(0)$ (see~\cite[Section 7]{eulergothic}). Its Euler characteristic is given by the arithmetic function $\tfrac{-1}{6c_{D}}\, e(D,6)$, that will be defined and studied in~\autoref{sec:numbertheory}.

With the square-discriminant case in mind, we will number the components according to the Hilbert modular surface where they live, that is the components $G_{d^{2}}^{r}$ and $\R[d^{2}]^{r}$ live in $X_{d^{2}}(\frakb_{r})$, where the possible ideals $\frakb_{r}$ are given by~\autoref{prop:Xcomponents}. This choice of representatives $r|6$ (cf.~\autoref{lem:ODmodules}) makes the formulae relating the discriminant $d^{2}$ of the Teichm\"{u}ller curve and the degree $m$ of the associated minimal torus coverings simpler (see~\autoref{thm:areadiscgothic}).

\

An obvious first remark is that, except in the genus 2 case, all these formulae are valid only for non-square discriminants. The reason for this is that the usual method to calculate the Euler characteristic of a Teichm\"{u}ller curve $C$ is the determination of its fundamental class $[\overline{C}]$ in some compactification of the Hilbert modular surface it lies on. In the square discriminant case, one needs to take into account some extra boundary divisors that do not appear in the non-square case. We will next show that the contribution of these extra components to the Euler characteristic of the Teichm\"{u}ller curve $G_{d^{2}}^{r}$ is negligible when computing asymptotics.

\begin{lemma}\label{lem:formulafornonsquare} Let $D=d^{2}$ be a square discriminant and $\frakb_{r}$ an ideal of norm 6 in $\cO_{d^{2}}$. There exists a constant $\epsilon$ depending only on $r$ and on $(d,6)$ such that 
\[ - \tfrac{3}{2} \,\chi(X_{d^{2}}(\frakb_{r})) - 2\,\chi(R_{d^{2}}^{r}) \le 
 \, \chi(G_{d^{2}}^{r}) \,
\le \bigl( - \tfrac{3}{2} + \tfrac{\epsilon}{d} \bigr)\,\chi(X_{d^{2}}(\frakb_{r})) - 2\,\chi(R_{d^{2}}^{r})\,.\]
\end{lemma}

\begin{proof}By~\cite[Thm. 8.1]{eulergothic}, the class of the divisor of the gothic modular form $\cG_{D}$ in $H^{2}(\widehat{X}_{D}(\frakb_{r});\bQ)$ agrees for non-square discriminant $D$ with
	\[ [\div(\cG_{D})] = [G_{D}(\frakb)] + 2[R_{D}(\frakb)]\,.\]

In the case of a square discriminant $D=d^{2}$ one just needs to determine the vanishing order of the modular form $\cG^{r}_{d^{2}}$ along the extra boundary components $S_{d^{2}}^{j}(\frakb_{r})$.

Let us first assume $\frakb_{r}=\frakb_{1}$, with associated $(1,6)$-symplectically adapted basis $\cO_{d^{2}}^{\vee} = \langle (\tfrac{1}{d},-\tfrac{1}{d}),(0,1)\rangle=\langle \eta_{1},\eta_{2}\rangle$ (see~\cite[\S 4]{eulergothic}). The corresponding linear forms read
	\[
	\rho(x_{1},x_{2}) = \frac{x_{1}}{d} \quad\mbox{ and }\quad
	\rho^{\sigma}(x_{1},x_{2}) = -\frac{x_{1}}{d} + x_{2}\,,
	\]
and one can write the following Fourier expansion of $\cG_{d^{2}}^{r}(\bm\tau)$ around the (unique!) cusp at infinity (cf.~\cite[Prop. 8.2]{eulergothic}, where the convention of a $(2,3)$-simplectically adapted basis was used )
	\begin{align*}
	\cG_{D}(\bm\tau)  = 8\pi^{2}i\cdot\Bigl( &\sum_{\substack{\bba\in\Lambda_{\frac{1}{2},0}  \\ \bbb\in\Lambda_{\frac{1}{2},\frac{1}{6}}}}\!\!\!  	
	(-1)^{a_{1}+b_{1}} \tfrac{(-a_{1}+da_{2})(-b_{1}+db_{2})}{d^{2}} \ 
	q_{1}^{\frac{a_{1}^{2}+b_{1}^{2}}{d^{2}}} 
	q_{2}^{(-\frac{a_{1}}{d}+a_{2})^{2}+(-\frac{b_{1}}{d}+b_{2})^{2} } 
	 \\
	- &\sum_{\substack{\bba\in\Lambda_{\frac{1}{2},\frac{1}{2}}  \\ \bbb\in\Lambda_{\frac{1}{2},\frac{2}{6}}}}\!\!\!  	
	(-1)^{a_{1}+b_{1}}  \tfrac{(-a_{1}+da_{2})(-b_{1}+db_{2})}{d^{2}} \ 
	q_{1}^{\frac{a_{1}^{2}+b_{1}^{2}}{d^{2}}} 
	q_{2}^{(-\frac{a_{1}}{d}+a_{2})^{2}+(-\frac{b_{1}}{d}+b_{2})^{2} } \Bigr)
	\,,
	\end{align*}
where $\Lambda_{i,j}\coloneqq \bZ^{2}+(i,j)$ and $q_{k}=\exp(\pi i \tau_{k})$.


The smallest $q_{1}$-exponent is $\tfrac{1}{2d^{2}}$, achieved by all terms $\bba=(\pm\tfrac{1}{2},a_{2})$, $\bbb=(\pm\tfrac{1}{2},b_{2})$. By looking at the smallest $q_{2}$-exponents of its coefficient
	\begin{align*}
	2\sum_{a,b\in\bZ} \biggl(
	\ \ &q_{2}^{(\frac{-1}{2d}+a)^{2}+(\frac{1}{2d}+b+\frac{1}{6})^{2}} \cdot
	\tfrac{(-\tfrac{1}{2}+da)(\tfrac{1}{2}+db+\tfrac{d}{6})}{d^{2}} \\[-0.3cm]
	- &q_{2}^{(\frac{-1}{2d}+a)^{2}+(\frac{-1}{2d}+b+\frac{1}{6})^{2}} \cdot
	\tfrac{(-\tfrac{1}{2}+da)(-\tfrac{1}{2}+db+\tfrac{d}{6})}{d^{2}} \\
	- &q_{2}^{(\frac{-1}{2d}+a+\tfrac{1}{2})^{2} + (\frac{1}{2d}+b+\frac{2}{6})^{2}} \cdot 
	\tfrac{(-\tfrac{1}{2}+da+\tfrac{d}{2})(\tfrac{1}{2}+db+\tfrac{2d}{6})}{d^{2}} \\
	+ &q_{2}^{(\frac{-1}{2d}+a+\tfrac{1}{2})^{2} + (\frac{-1}{2d}+b+\frac{2}{6})^{2}} \cdot 
	\tfrac{(-\tfrac{1}{2}+da+\tfrac{d}{2})(-\tfrac{1}{2}+db+\tfrac{2d}{6})}{d^{2}}
	\biggr)
	\end{align*}
it is easy to see that this coefficient does not vanish generically, and therefore the modular form behaves along $S_{d^{2}}^{2}(\frakb_{1})$ as $q_{1}^{1/(2d^{2})}$. 

Now, the stabiliser around $(\infty,\infty)$ is given by the matrix group
	\[\Stab_{(\infty,\infty)}(\SL(\frakb_{1}\oplus\cO_{d^{2}}^{\vee}))=\left\{ \begin{pmatrix} 1 & \mu \\ 0 & 1 \end{pmatrix}\,:\, \mu\in M = \sqrt{D}\frakb_{1}\right\}\,.\]
Consider the basis $M^{\vee}= \tfrac{1}{6d^{2}}\frakb_{6} = \tfrac{1}{6d^{2}}\langle (6,6),(0,d)\rangle\eqqcolon \langle (\alpha_{1},\alpha_{2}),(\beta_{1},\beta_{2})\rangle$. The exponentials
	\[\left\{\begin{array}{l}
	X = q_{1}^{2\alpha_{1}}q_{2}^{2\alpha_{2}} = q_{1}^{2/d^{2}}q_{2}^{2/d^{2}}\,,\\
	Y = q_{1}^{2\beta_{1}}q_{2}^{2\beta_{2}} = q_{2}^{1/(3d)}
	\end{array}\right.\]
are local parameters around $(\infty,\infty)$ for the Hilbert modular surface $X_{d^{2}}(\frakb_{1})$. Replacing $(q_{1},q_{2})$ by the local parameters $(X=q_{1}^{2/d^{2}}q_{2}^{2/d^{2}},Y=q_{2}^{1/(3d)})$ the correct minimal exponent becomes $X^{1/4}$, and the modular form vanishes to order $1/4$ along $S_{d^{2}}^{2}(\frakb_{1})$.

Finally, the Euler characteristic of $G_{d^{2}}^{1}$ can be written as
	\[\chi(G_{d^{2}}^{1}) = -[\omega_{1}]\cdot [G_{d_{2}}^{1}] = -[\omega_{1}]\cdot\left( \frac{1}{2}[\omega_{1}]+\frac{3}{2}[\omega_{2}]-2[R_{d^{2}}^{1}] - \frac{1}{4}[S_{d^{2}}^{2}(\frakb_{1})] \right) \]

Pairing $-[\omega_{1}]$ with $[S_{d^{2}}^{2}(\frakb_{1})]$ computes the Euler characteristic of $S_{d^{2}}^{2}(\frakb_{1})$, which by~\autoref{lem:componentsboundary} is bounded by (see for example~\cite[Prop. 10.5]{bainbridge07})
	\[\chi(S_{d^{2}}^{2}(\frakb_{1})) \le 12\cdot \chi(\Gamma_{1}(2160\, d)\backslash \bH) = - 2160^{2} d^{2} \sum_{r|30\,d}\frac{\mu(r)}{r^{2}}\,.\]
By~\eqref{eq:chiXd2} this is in turn bounded by $\frac{\epsilon}{d}X_{d^{2}}(\frakb_{1})$ for some constant $\epsilon$.

The other cases follow from doing the same calculations for each of the components of $S_{d^{2}}^{2}(\frakb_{r})$.
\end{proof}


\begin{remark}This approach would give a precise formula for the Euler characteristics of arithmetic gothic Teichm\"{u}ller curves. However, the case by case analysis of the different components of $S_{d^{2}}^{2}(\frakb_{r})$ for each of the possible values of $r$ and of $(d,6)$ and the corresponding calculation of the vanishing order of $G_{d^{2}}^{r}$ exceed the objectives of this paper. For example, computer simulations show that for $r=1$ the Euler characteristic of $G_{d^{2}}^{1}$ is given by:
	\begin{align*}
	\chi(G_{d^{2}}^{1})
	&=
	\left\{\begin{array}{ll}
	(\frac{2}{d}-\frac{3}{2})\chi(X_{d^{2}}(\frakb_{1})) -2\chi(R_{d^{2}}^{1})\,, & \mbox{if $(6,d)=1$,} \\
	(\frac{6}{d}-\frac{3}{2})\chi(X_{d^{2}}(\frakb_{1})) -2\chi(R_{d^{2}}^{1})\,, & \mbox{if $(6,d)=2$,} \\
	(\frac{3}{d}-\frac{3}{2})\chi(X_{d^{2}}(\frakb_{1})) -2\chi(R_{d^{2}}^{1})\,, & \mbox{if $(6,d)=3$,} \\
	(\frac{9}{d}-\frac{3}{2})\chi(X_{d^{2}}(\frakb_{1})) -2\chi(R_{d^{2}}^{1})\,, & \mbox{if $(6,d)=6$.} \\
	\end{array}\right.
	\end{align*}
\end{remark}

\subsection{Minimal torus covers in Teichm\"{u}ller curves}

Our strategy is based on the fact that Euler characteristics allow us to compute the number of minimal torus covers on a given Teichm\"{u}ller curve.

\begin{lemma}\label{lem:stseulerchar} Let $\pi:(X,\omega)\to (\bC/Per(X,\omega),dz)$ be a minimal torus cover of degree $m$ and denote by $\widetilde{C}=\GL^{+}(2,\bR)(X,\omega)$ its orbit and by $C$ the corresponding Teichm\"{u}ller curve, that the projection of $\widetilde{C}$ to the moduli space. The number $\cS_{m,m}(\widetilde{C})$ of minimal torus covers of degree~$m$ and area~$m$ in $\widetilde{C}$ equals $-6\chi(C)$.
\end{lemma}

\begin{proof}By considering the flat surface $M^{-1}\cdot(X,\omega)\in \widetilde{C}$, where $M$ is the Hermite normal form of the lattice $\Per(X,\omega)$, we can assume that both the degree of $(X,\omega)\to (\bC/\Per(X,\omega),dz)$ and the area of $(X,\omega)$ are equal to $m$. In particular, since the action of a matrix $A\in \GL^{+}(2,\bR)$ changes the area by $\det(A)$, we can restrict our attention to its orbit under $\SL(2,\bR)$.

Now, each  point $A\cdot(X,\omega)$ with $A\in\SL(2,\bR)$ comes with an associated covering $A\cdot(X,\omega)\to A\cdot(E,dz)$ and one can therefore define a map $\overline{\pi}: \SL(2,\bR)\cdot (X,\omega) \to M_{1,1}$, which sends each point $A\cdot(X,\omega)$ to the elliptic curve $A\cdot (E,dz)$ it covers. The number of minimal torus covers on $\widetilde{C}$ of degree~$m$ and area~$m$ is precisely the number of preimages of the square torus $E$. Since this covering is unramified, this is given by the degree of the cover $\overline{\pi}$ which agrees with $\chi(C)/\chi(M_{1,1}) = -6 \chi(C)$.
\end{proof}

We can now prove the theorem at the beginning of the section. 

\begin{proof}[Proof of~\autoref{thm:areadiscgothic}]Let $\pi_{min}:(X,\omega)\to (\bC/\Per(X,\omega),dz)$ be a minimal torus cover in $\G$ of degree and area~$m$, and write $E_{X}=\bC/\Per(X,\omega)$. Accordingly, let $\Prym X$ correspond to a point $\bm\tau\in\bH^{2}$ in the Hilbert modular variety, that is $T_{\bm\tau}=\Prym X \in X_{d^{2}}(\frakb_{r})$ for some $d$ and $r$.

By the general theory, the covering map $\pi_{min}$ induces an inclusion $\pi_{min}^{*}: E_{X}\to \Jac X$, where the degree of the cover $m=\deg(\pi_{min})$ agrees with the type of the restriction $\cO_{\Jac X}(\Theta)|_{\pi_{min}^{*}E_{X}}=m\cO_{E_{X}}(0)$ of the principal polarisation to $\pi_{min}^{*}E_{X}$. 

With the usual normalisation, we always assume that the eigenform $\omega=\pi_{min}^{*}dz$ corresponds to $du_{1}$ in $X_{d^{2}}(\frakb_{r})$, and therefore $\pi_{min}^{*}E_{X}=E_{1}$ as in~\autoref{prop:ellipticcurves} (one can equivalently see the Teichm\"{u}ller curve as a curve in $X_{d^{2}}(\frakb_{r}^{\sigma})$ with $du_{2}$ as an eigenform). The theorem implies that, for the component $G_{d^{2}}^{r}$ belonging to the Hilbert modular surface $X_{d^{2}}(\frakb_{r})$, the degree is given by $m=\lcm(d,r)$ (see~\autoref{tab:areadiscriminant} for the different values appearing in each case). Writing $d$ in terms of $m$ yields the result. 
%
%

\end{proof}

\begin{table}[h]\centering
\begin{tabular}[t]{ll}
	\begin{tabular}[t]{| l | c | l |} 
	\hline 
	$\bm{d\bmod{6}}$ & $\bm{\frakb_{r}}$ & $\quad\bm{n}$ \\ [0.5ex]  \hline \hline 
	$0$ & $\frakb_{1}$ & $n=d$ \\ \hline
	$3$ & $\frakb_{1}$ & $n=d$ \\
	    & $\frakb_{2}$ & $n=2d$ \\ \hline
	$2,4$ & $\frakb_{1}$ & $n=d$ \\ 
	      & $\frakb_{3}$ & $n=3d$ \\ \hline
	$1,5$ & $\frakb_{1}$ & $n=d$ \\ 
	      & $\frakb_{2}$ & $n=2d$ \\
	      & $\frakb_{3}$ & $n=3d$ \\
	      & $\frakb_{6}$ & $n=6d$ \\ \hline
	\end{tabular}
&
	\begin{tabular}[t]{| l | c | l |} 
	\hline 
	$\bm{d\bmod{2}}$ & $\bm{\frakb_{r}}$ & $\quad\bm{n}$ \\ [0.5ex]  \hline \hline 
	$0$ & $\frakb_{1}$ & $n=d$ \\ \hline
	$1$ & $\frakb_{1}$ & $n=d$ \\
	    & $\frakb_{2}$ & $n=2d$ \\ \hline
	\end{tabular}
\\
&
\end{tabular}
\caption{Relation between the degree $n$ of a minimal torus cover and the discriminant $d^{2}$ of the Teichm\"{u}ller curve for $\G$ and $\cP_{4}$.}\label{tab:areadiscriminant}
\end{table}

The very same strategy proves also the results of McMullen and Lanneau-Nguyen in genus 2, 3 and 4.

\begin{theorem}[\cite{mcmTCingenustwo,lanneaunguyen,LNPrym4}]\label{thm:areadiscprym} For each $m>0$, the sets of minimal torus covers of degree and area $m$ in the stratum $\cH(2)$ and in the Prym loci $\cP_{3}$ and $\cP_{4}$ distribute among different Teichm\"{u}ller curves in the following way:
	\begin{align*}
	\bullet&\ \cS_{m,m}(\cH(2))  && \hspace{-0.5cm} = \cS_{m,m}(W_{m^{2}}(2))\,; \\[4pt]
	\bullet&\ \cS_{m,m}(\cP_{3}) && \hspace{-0.5cm} = \left\{
	\begin{array}{ll}
	\cS_{m,m}(W^{1}_{m^{2}}(4))\,, & \mbox{if $m\not\equiv 2 \bmod{4}$,} \\
	\cS_{m,m}(W^{1}_{m^{2}}(4)) \cup \cS_{m,m}(W^{2}_{(m/2)^{2}}(4))\,,  &\mbox{if $m\equiv 2\bmod{4}$;} 
	\end{array}\right.\\[4pt]
	\bullet&\ \cS_{m,m}(\cP_{4}) && \hspace{-0.5cm} = \left\{
	\begin{array}{ll}
	\cS_{m,m}(W_{(m/2)^{2}}(6))\,, &\mbox{if $m\equiv 0\bmod{2}$,} \\
	\varnothing\,,  &\mbox{if $m\equiv 1\bmod{2}$.}
	\end{array}\right.
	\end{align*}

\end{theorem}

%% file: sec_number.tex
\section{Asymptotics of sums over arithmetic progressions} \label{sec:numbertheory}

In this section we use Dirichlet series and modular forms to calculate the asymptotic behaviour of certain functions closely related to our counting problem. More precisely, the volumes of the affine invariant manifolds in which we are interested will be given by an infinite sum of Euler characteristics of Teichm\"{u}ller curves weighted by certain $\sigma$ function. These can in turn be written in terms of the Euler characteristics of Hilbert modular surfaces and of the reducible locus, which are given by some arithmetic functions $e(D,k)$.

However, the formulae appearing in the gothic case depend heavily on the congruence class of the discriminant considered (see~\autoref{thm:areadiscgothic} and~\autoref{thm:eulergothic}). This implies that the sum for the volume will split into different sums running through different arithmetic progressions, and we will need to estimate each of the summands separately.

For a quadratic discriminant $D=f^{2}D_{0}$ with conductor $f$ and a positive integer $k\ge 1$, let us define the set of prototypes
\begin{align*}
\cP_k(D) \=  \bigl\{ &[a,b,c] \in \bZ^{3}\, :\,
\mbox{$a>0>c$ , $D=b^{2} - 4\cdot k \cdot ac$} \\
& \mbox{and $\gcd\left(f,b,c_{0}\right)=1$, where $c=c_{0}^{2}\cdot c'$ and $c'$ is square-free }\bigr\}\,.\nonumber
\end{align*} 

The Euler characteristics of $X_{D}(\frakb)$ and $\R(\frakb)$ are intimately related to the following arithmetic functions
\bas
e(1,k) &= -\frac{1}{12}\,, \\
e(D,k) &= \!\!\! \sum_{[a,b,c]\in \cP_{k}(D)} \!\!\! a\,, \\
\a(d) &= |\SL_{2}(\mathbb{Z}/d\mathbb{Z})| = d\sum_{m|d}\mu\bigl(\tfrac{d}{m}\bigr)m^{2}\,.
\eas

The asymptotics of $e(D,k)$ for fundamental discriminants $D$ were calculated in~\cite[Thm. 10.1]{eulergothic}. Our main objective is to extend this study to the case $D=d^{2}$ and to estimate the growth of some related convolutions.

The main results of this section are the following propositions.

\begin{prop}\label{prop:e(d,k)} The asymptotic behaviour of the functions $e(d^{2},1)$ and $e(d^{2},6)$ as $d$ grows is given by:
	\bas
	e(d^{2},1) &= \frac{5}{12}\,\a(d) +  O(d^{5/2}) \,, \\
	e(d^{2},6) &= 	
	\begin{cases}
	\hfill 2\cdot \tfrac{1}{60} \cdot \a(d) +  O(d^{5/2})\,,  & \text{if $(6,d)=1$,} \\
	\tfrac{3}{2}\cdot \tfrac{1}{60} \cdot \a(d) +  O(d^{5/2})\,,  & \text{if $(6,d)=2$,} \\
	\tfrac{4}{3}\cdot \tfrac{1}{60} \cdot \a(d) +  O(d^{5/2})\,,  & \text{if $(6,d)=3$,} \\
	\hfill \tfrac{1}{60} \cdot \a(d) +  O(d^{5/2})\,,  & \text{if $(6,d)=6$.}
	\end{cases}
	\eas
\end{prop}

As stated above, the summands appearing in the formulae for the Masur-Veech volumes will be formed in our case by certain convolutions where the divisors run over certain congruence class. Let us define the following sums
	\be \label{eq:sum_convolutions}
	\SS_{k}(D) = \sum_{d=1}^{D} \!\!\! \sum_{ \genfrac{}{}{0pt}{}{m|d}{(m,k)=k} } \!\!\! \sigma\bigl(\tfrac{d}{m}\bigr) \a(m) =
	\sum_{d=1}^{D} \!\!\! \sum_{ \genfrac{}{}{0pt}{}{m|d}{(m,k)=k} } \!\!\! \sum_{r|m} \tfrac{m^{3}}{r^{2}}\sigma\bigl(\tfrac{d}{m}\bigr)\mu(r)\,,
	\ee
We note here that the usual convolution $(\sigma\ast \a)(n)$ is simply $\sigma_{3}(n)$. In our case it will be enough to study the behaviour of $S_{k}(D)$ when $k$ is~1, a prime or a product of two primes.

\begin{prop}\label{prop:S_k} The asymptotic behaviour of the sums $S_{k}(D)$ as $D$ grows is given in the cases $k=1,p,pq$ by:
	\bas
	\SS_{1}(D)  &= \frac{\pi^{4}D^{4}}{2^{3}\cdot 3^{2}\cdot 5} + O(D^{3})\,, \\
	\SS_{p}(D)  &= \frac{\pi^{4}D^{4}}{2^{3}\cdot 3^{2}\cdot 5}\cdot\frac{p+1}{p^{2}+p+1} + O(D^{3})\,, \\
	\SS_{pq}(D) &= \frac{\pi^{4}D^{4}}{2^{3}\cdot 3^{2}\cdot 5}\cdot \frac{p+1}{p^{2}+p+1} \cdot \frac{q+1}{q^{2}+q+1} + O(D^{3})\,,
	\eas
for any different primes $p$ and $q$.
\end{prop}

\subsection{Dirichlet series}

Given an arithmetic function $f$, its Dirichlet series is the formal series given by
	\[D_{f}(s)=\sum_{n=1}^{\infty} \frac{f(n)}{n^{s}}\,,\quad \mbox{where $s\in\mathbb{C}$}.\]
Although we are in general not interested in the convergence of the series, it is worth noting that if two Dirichlet series $D_{f}(s)$ and $D_{g}(s)$ agree on some half-plane $\{\Re(s)>a\}$ where they converge absolutely, then $f(n)=g(n)$ for all $n\in\mathbb{N}$.

Dirichlet series are useful to calculate convolutions of arithmetic functions. More precisely, for multiplicative functions $f$ and $g$, we define its convolution as the sum
	\[(f * g)(d) = \sum_{m n =d} f(m)g(n)\,.\]
Then one has
	\[D_{f * g}(s) = D_{f}(s) D_{g}(s)\,.\]

Convolution is commutative and associative. In particular, the equality of Dirichlet series above can be generalised to an arbitrary number of factors in the convolution. 

Note also that, if we denote by $j_{k}(n)=n^{k}$ the $k$-th power and by $\bm{1}(n)=1$ the constant function 1, then $D_{j_{k}\cdot f}(s) = D_{f}(s-k)$ and $D_{\bm{1}}(s) = \zeta(s)$.

The following equalities are either well known or straightforward (see for example~\cite[\S 11.5]{Apo}):
	\begin{align}
	\begin{split}
	D_{\sigma_{k}}(s) &= \zeta(s)\zeta(s-k)\,,\quad \mbox{for $k>1$}\,, \\ \label{eq:dirichletseries}
	D_{\mu}(s) &= \zeta(s)^{-1}\,, \\
	D_{j_{k}}(s) &= \zeta(s-k)\,, \\
	D_{j\cdot\mu}(s) &= \zeta(s-1)^{-1} \,.
	\end{split}
	\end{align}

In particular, as noted above, one has
	\be \label{eq:convolution_sigma_and_a}
	D_{\a}(s)=\frac{\zeta(s-3)}{\zeta(s-1)} \quad\mbox{and}\quad 
	D_{\a\ast \sigma}(s)=\zeta(s)\zeta(s-3) = D_{\sigma_{3}}(s)\,.
	\ee

\subsection{Divisor sums and modular forms}

%
%
%

Our analysis is based on the results of Zagier in \cite[Section~4]{DonNZeta}. Recall that the theta series $\theta$ and the Eisenstein series $G_{2}$ are modular forms defined by
\bes
\theta(\tau) \= \sum_{\ell = -\infty}^\infty e^{\pi i \ell^2 \tau},
\quad
G_{2}(\tau) \= \frac{-1}{24} \,+\,
\sum_{a=1}^\infty \sigma(a) e^{2\pi i a\tau} \,.
\ees
Then the coefficients of the modular form 
\bes
F_{k}(\tau) \,:=\, G_2(2k\tau) \theta(\tau)
\= \sum_{n=0}^\infty e_{k}(n) e^{\pi i n\tau} \,.
\ees
are given by
	\[e_{k}(n)= \sum_{\genfrac{}{}{0pt}{}{b^2 \equiv n\mod 4k,}{|b| \leq \sqrt{n}}}
\sigma\Bigl(\frac{n-b^2}{4k}\Bigr)\,,\]
where we define $\sigma(0)=-\tfrac{1}{24}$. It is immediate to see that, for a discriminant $D=f^{2}D_{0}$ as above, one has
	\begin{equation} \label{eq:e_and_a}
	e_{k}(D) = \sum_{m|f} e(\tfrac{D}{m^{2}},k)\,.
	\end{equation}

One can therefore use the coefficients of this modular form to determine the values $e(D,k)$ using M\"{o}bius inversion (see~\autoref{prop:e(d,k)} below). This will allow us to study the asymptotic behaviour of $e(d^{2},1)$ and $e(d^{2},6)$, and therefore of the Euler characteristics in which we are interested. 

Following Zagier, let us introduce the following Gauss sum (see~\cite[Thm.2]{DonNZeta} for the definitions and the facts claimed below)
\bes
\gamma_c(n) = c^{-1/2} \sum_{a=1}^{2c} \lambda(a,c) e^{-\pi i na/c}
\ees
where $\lambda(a,c)$ is a Legendre symbol times a power of~$i$, depending on the parities of~$a$ and~$c$. This is a multiplicative function in $c$, given at prime powers $c=p^{r}$ for $n=d^{2}$ a square by:
\be \label{eq:gamma2}
\gamma_{2^r}(d^{2}) = 
\begin{cases}
1 & \text{if $r=0$,} \\
2^{r/2} & \text{if $r$ is even and $\nu_{2}(d^{2})=r-2$,} \\
2^{(r-1)/2} & \text{if $r$ is odd and $\nu_{2}(d^{2})\ge r-1$,} \\
0 & \text{otherwise,} \\
\end{cases} 
\ee
for $p=2$ and
\be \label{eq:gammap}
\gamma_{p^r}(d^{2}) = 
\begin{cases}
1 & \text{if $r=0$,} \\ 
p^{r/2-1}(p-1) & \text{if $r$ is even and $\nu_{p}(d^{2})\ge r$,} \\
p^{(r-1)/2} & \text{if $r$ is odd and $\nu_{p}(d^{2})= r-1$,} \\
0 & \text{otherwise,} \\
\end{cases} 
\ee
for odd primes, where $\nu_{p}(m)$ denotes the $p$-adic valuation of $m$. These Gauss sums allow us to define the coefficients 
\ba \label{eq:defstarenk}
e^{*}_{k}(n) &= \sum_{c=1}^\infty \frac{(c,2k)^2}{c^{2}}\gamma_c(n)\,, \\
\overline{e}_{k}(n) &=\frac{\pi^{2}}{72 k^{2}}n^{3/2} e^{*}_{k}(n)\,,
\ea
which are well defined since the summands grow at most like $c^{-2}$ (this is clear for $n=d^{2}$, for the general case it follows from the expressions in~\cite[Thm. 2]{DonNZeta}).

These coefficients were introduced by Zagier in~\cite[\S 4]{DonNZeta} for the case $k=1$. There he defines the Dirichlet series $E_{n}(s)$, which for square $n=d^{2}$ has the form
	\be \label{eq:DirichletSeries}
	E_{d^{2}}(s) = \frac{\zeta(s)}{\zeta(2s)} \sum_{ \genfrac{}{}{0pt}{}{a,c\ge 1}{ac|d} }\frac{\mu(a)}{c^{2s-1}a^{s}}
	\ee
and whose value at $s=2$ is related to the coefficients above by $\overline{e}_{1}(n)=\frac{\pi^{2}}{36}n^{3/2}E_{n}(2)$. It can also be proved (see~\cite{DonNZeta} for the case $k=1$ and the proof of Theorem 10.1 in~\cite{eulergothic} and the references therein for the general case) that
	\be\label{eq:e_k_approx}
	e_{k}(n)=\overline{e}_{k}(n)+ O(n^{5/4})\quad \mbox{ as $n\to \infty$.}
	\ee

All this allows us to prove the following result.
\par
\begin{lemma} \label{lem:e1_and_ek}
Let $k$ be a square-free positive integer. Then 
\bes
e^{*}_{k}(d^{2}) =  k^{2} \sum_{m|k} \mu(m)\Bigl(\!\! \prod_{\genfrac{}{}{0pt}{}{p | m}{\text{$p$ prime}}}\!\! \tfrac{p^{2}-1}{p^{2}+1} \Bigr)\,e^{*}_{1}(d_{m}^{2})  \,,
\ees
where 
	\[ d_{m}\coloneqq \frac{d}{\prod_{p|m} p^{\nu_{p}(d)}} = \max\{x|d\,:\, (x,m)=1\}\,.\]
\end{lemma}

\begin{proof} The summands in the function $e_{k}^{*}(n)$ are weakly multiplicative in $c$, so it admits an Euler product expansion
	\[e_{k}^{*}(n)= \prod_{p\ \text{prime}} P_{k}(p,n)\,,\quad \mbox{where}\quad 
	P_{k}(p,n)=1 + \sum_{j=1}^{\infty} \frac{(p^{j},2k)^{2}}{p^{2j}}\gamma_{p^{j}}(n)\,.\]
In particular, $P_{k}(p,n)=P_{1}(p,n)$ whenever $p\nmid k$, whereas for $p\mid k$ one has
	\bas
	P_{k}(p,n) &= p^{2}\, P_{1}(p,n) - (p^{2}-1)\,,\quad \mbox{if $p\neq 2$ and} \\
	P_{k}(2,n) &= 4\, P_{1}(2,n) - 3 - 3\gamma_{2}(n)\,.
	\eas

Let us focus on the case $n=d^{2}$ a square. If $p\nmid k$ and $p\neq 2$, Equation~\eqref{eq:gammap} together with the fact that $p\nmid d_{p}$ imply that $P_{1}(p, d^{2}_{p})= 1 +1/p^{2}$, and therefore 
	\[P_{k}(p,d^{2}) = p^{2} \Bigl( P_{1}(p,d^{2}) - \tfrac{(p^{2}-1)}{(p^{2}+1)}\, P_{1}(p, d^{2}_{p})\Bigr)\,.\]

As for $p=2$, applying Equation~\eqref{eq:gamma2} one has $\gamma_{2}(d^{2})=1$ and $P_{1}(2,d_{2}^{2})= 5/2$ and hence, if $2\nmid k$,
	\[P_{k}(2,d^{2}) = 4 \Bigl( P_{1}(2,d^{2}) - \tfrac{3}{5}\, P_{1}(2, d^{2}_{2})\Bigr)\,.\]

Finally, we can write
	\bas
	e^{*}_{k}(d^{2}) &= \prod_{p \nmid k} P_{1}(p,d^{2}) \cdot \prod_{p \mid k} p^{2} \Bigl( P_{1}(p,d^{2}) - \tfrac{p^{2}-1}{p^{2}+1}\, P_{1}(p, d^{2}_{p})\Bigr) \\
	&= k^{2}\sum_{m\mid k} \Bigl(  \mu(m) \prod_{p \nmid m} P_{1}(p,d^{2}) \cdot \prod_{p|m} \tfrac{p^{2}-1}{p^{2}+1}\, P_{1}(p, d^{2}_{p})\Bigr) \,,
	\eas
which gives the desired formula since, for $p\mid m$, the equality $P_{1}(p, d^{2}_{p})=P_{1}(p, d^{2}_{m})$ holds.
\end{proof}

In the particular case that we are interested the previous lemma yields the following result.

\begin{cor} \label{cor:e1_and_e6}
\bes
e^{*}_{6}(d^{2}) = 36\Bigl( e^{*}_{1}(d^{2})  - \tfrac{3}{5} \, e^{*}_{1}(d_{2}^{2}) - \tfrac{4}{5} \, e^{*}_{1}(d_{3}^{2})  + \tfrac{12}{25} \, e^{*}_{1}(d_{6}^{2})  \Bigr)\,.
\ees
\hfill$\square$
\end{cor}

\subsection{Asymptotic behaviour of the Euler characteristics}

The Euler characteristics of $X_{D}(\frakb)$ and $\R(\frakb)$ can be written in terms of the arithmetic functions defined above. By~\cite[\S 1]{Hir} and~\cite[Lemma 7.5]{eulergothic}, they are given by:
	\begin{align}
	\begin{split}
	\chi(X_{D}) &= \tfrac{1}{30}\, e(D,1)\,,\quad\mbox{if $D$ is not a square,} \label{eq:eulercharande} \\
	\chi(X_{d^{2}}) &= \tfrac{1}{72}\, \a(d)\,, \\
	\chi(R_{D}^{r}) &= \tfrac{-1}{6c_{D}}\, e(D,6)\,, \\
	\end{split}
	\end{align}
where $c_{D}$ denotes the number of ideals of norm $6$ in $\cO_{D}\otimes \bQ$. In the case of a square discriminant $D=d^{2}$ one has $c_{d^{2}}=\sigma_{0}\bigl(\frac{6}{(d,6)}\bigr)$. 

We can finally use the results of the previous section to estimate the asymptotics of the functions $e(d^{2},k)$ and $S_{k}(D)$.


\begin{proof}[Proof of~\autoref{prop:e(d,k)}] Applying the formula for $\overline{e}_{1}(n)$ in terms of $E_{n}(2)$ to the particular case of $n=d^{2}$, one has 
	\[\overline{e}_{1}(d^{2})=\frac{5}{12}d^{3} \sum_{ \genfrac{}{}{0pt}{}{a,c\ge 1}{ac|d} }\frac{\mu(a)}{c^{3}a^{2}}\,.\]
Now M\"{o}bius inversion applied to~\eqref{eq:e_and_a} yields
	\bas
	e(d^{2},1) &= \sum_{m|d} \mu\bigl(\tfrac{d}{m}\bigr) e_{1}(m^{2}) = \sum_{m|d} \mu\bigl(\tfrac{d}{m}\bigr) \overline{e}_{1}(m^{2}) + \sum_{m|d} \mu\bigl(\tfrac{d}{m}\bigr) (e_{1}(m^{2}) - \overline{e}_{1}(m^{2})) = \\
	&= \frac{5}{12} \sum_{m|d} \mu\bigl(\tfrac{d}{m}\bigr) m^{3} \sum_{ \genfrac{}{}{0pt}{}{a,c\ge 1}{ac|m} }\frac{\mu(a)}{c^{3}a^{2}} + \sum_{m|d} \mu\bigl(\tfrac{d}{m}\bigr) (e_{1}(m^{2}) - \overline{e}_{1}(m^{2}))\,.
	\eas

The first summand is given by the quadruple convolution
	\[\frac{5}{12} \sum_{m|d} \mu\bigl(\tfrac{d}{m}\bigr) m^{3} \sum_{ \genfrac{}{}{0pt}{}{a,c\ge 1}{ac|m} }\frac{\mu(a)}{c^{3}a^{2}} = 
	\frac{5}{12} \sum_{acrs=d} \mu(s) \mu(a)a r^{3} = \tfrac{5}{12} \a(d)\,, \]
where we have used the fact that the Dirichlet series of the quadruple convolution agrees with $D_{a(n)}(s) = \zeta(s-3)/\zeta(s-1)$.

As for the second summand, M\"{o}bius inversion together with the estimate~\eqref{eq:e_k_approx} show that it is $O(d^{5/2})$ as $d$ grows.

\

The same reasoning applied to the case $k=6$ yields
	\[
	e(d^{2},6) 
	= \sum_{m|d} \mu\bigl(\tfrac{d}{m}\bigr) e_{6}(m^{2}) 
	=\sum_{m|d} \mu\bigl(\tfrac{d}{m}\bigr) \overline{e}_{6}(m^{2})  + \sum_{m|d} \mu\bigl(\tfrac{d}{m}\bigr) (e_{6}(m^{2})- \overline{e}_{6}(m^{2}))\,.\]

Now, by~\autoref{cor:e1_and_e6} the first summand equals
	\[
	\sum_{m|d} \mu\bigl(\tfrac{d}{m}\bigr) \overline{e}_{6}(m^{2}) = 
	\tfrac{\pi^{2}}{72} \sum_{m|d} \mu\bigl(\tfrac{d}{m}\bigr)m^{3} \Bigl( e^{*}_{1}(m^{2})  - \tfrac{3}{5} \, e^{*}_{1}(m_{2}^{2}) - \tfrac{4}{5} \, e^{*}_{1}(m_{3}^{2})  + \tfrac{12}{25} \, e^{*}_{1}(m_{6}^{2})  \Bigr)\,,\]

The result follows by a direct application of the two lemmas 
below to each of the four summands. As an example we next show the case $(6,d)=2$.

By applying again M\"{o}bius inversion to~\eqref{eq:e_and_a} and using~\autoref{lem:technicallemma} below, these four summands can be written as
	\bas
	\tfrac{\pi^{2}}{72} \sum_{m|d} \mu\bigl(\tfrac{d}{m}\bigr)m^{3} e^{*}_{1}(m^{2}) &= e(d^{2},1) +  O(d^{5/2})\,, \\
	\tfrac{3}{5} \tfrac{\pi^{2}}{72} \sum_{m|d} \mu\bigl(\tfrac{d}{m}\bigr)m^{3} e^{*}_{1}(m_{2}^{2}) &= 
	\tfrac{21}{5} \cdot 2^{3\nu_{2}(d)-3} \cdot e(d_{2}^{2},1) +  O(d^{5/2})\,,
	\\
	\tfrac{4}{5} \tfrac{\pi^{2}}{72} \sum_{m|d} \mu\bigl(\tfrac{d}{m}\bigr)m^{3} e^{*}_{1}(m_{3}^{2}) &= 
	\tfrac{4}{5} \cdot e(d^{2},1) +  O(d^{5/2})\,,
	\\
	\tfrac{12}{25} \tfrac{\pi^{2}}{72} \sum_{m|d} \mu\bigl(\tfrac{d}{m}\bigr)m^{3} e^{*}_{1}(m_{6}^{2}) &= 
	\tfrac{84}{25} \cdot 2^{3\nu_{2}(d)-3} \cdot e(d_{2}^{2},1) +  O(d^{5/2}) \,.
	\\
	\eas

We can now use the equality $e(d^{2},1) = 5/12\cdot \a(d) +  O(d^{5/2}) $ and~\autoref{lem:a(d)_and_a(d_p)} to get
	\bas
	\sum_{m|d} \mu\bigl(\tfrac{d}{m}\bigr) \overline{e}_{6}(m^{2}) &= 
	\bigl(\tfrac{5}{12} - \tfrac{1}{3}\bigr)\, \a(d) + 
	\bigl(-\tfrac{7}{4} + \tfrac{7}{5} \bigr)  2^{3\nu_{2}(d)-3} \cdot \a(d_{2}) +  O(d^{5/2}) = \\
	&= \tfrac{1}{12}\, \a(d) - \tfrac{7}{20} \cdot \tfrac{1}{6} \cdot \a(d)  +  O(d^{5/2})  = \tfrac{1}{40}\, \a(d)  +  O(d^{5/2}) \,.
	\eas
	
\end{proof}

In particular, the growth rates of the functions determining the Euler characteristics in~\eqref{eq:eulercharande} are the same in the square and non-square case.

\begin{proof}[Proof of~\autoref{prop:S_k}]The case $k=1$ is straightforward: by~\eqref{eq:convolution_sigma_and_a} and using the hyperbola trick one gets
	\bas
	S_{1}(D) &= \sum_{d=1}^{D} (\sigma\ast\a)(d) = \sum_{d=1}^{D} \sigma_{3}(d) =
	\sum_{x=1}^{D}\sum_{y=1}^{\lfloor\frac{D}{x}\rfloor} y^{3} = \sum_{x=1}^{D} \frac{1}{4} \left\lfloor\frac{D}{x}\right\rfloor^{2} \left(\left\lfloor\frac{D}{x}\right\rfloor+1\right)^{2}\,.
	\eas

Since the sum of the first $N$ cubes is given by $\frac{1}{4}N^{2}(N+1)^{2}$ and $\lfloor z\rfloor = z +O(1)$, this sum is dominated by $\frac{D^{4}}{4}\sum_{x\le D}\frac{1}{x^{4}}$, with an error term $O(D^{3})$. Since the sum of reciprocals of fourth powers tends to $\pi^{4}/90$, the result follows.

Now, for a fixed $D$ one has by~\autoref{lem:a(d)_and_a(d_p)}
	\bas
	S_{p}(D) &= \sum_{d=1}^{\lfloor D/p \rfloor}\sum_{m|d}\sigma\bigl(\tfrac{d}{m}\bigr) \a(pm) =\\
	&= p(p^{2}-1)\sum_{d=1}^{\lfloor D/p \rfloor}\sum_{m|d}\sigma\bigl(\tfrac{d}{m}\bigr) \a(m) + p \sum_{d=1}^{\lfloor D/p \rfloor} \!\!\!\! \sum_{ \genfrac{}{}{0pt}{}{m|d}{m\equiv 0 \bmod{p}} }\!\!\!\!\!  \sigma\bigl(\tfrac{d}{m}\bigr) \a(m) = \\
	&= p(p^{2}-1)\, S_{1}\bigl( \lfloor\tfrac{D}{p}\rfloor \bigr) + p\, S_{p} \bigl( \lfloor\tfrac{D}{p}\rfloor \bigr)\,.
	\eas
Since the growth of $S_{p}(D)$ is dominated by that of $S_{1}(D)= D^{4}\pi^{4}/360+ O(D^{3})$, the result follows from estimating $S_{p}(D)$ and $S_{p}(\lfloor D/p\rfloor)$.

Let now $p$ and $q$ be different primes. The same strategy as in the previous case yields
	\bas
	S_{pq}(D) =& \sum_{d=1}^{\lfloor D/pq \rfloor}\sum_{m|d}\sigma\bigl(\tfrac{d}{m}\bigr) \a(pqm) = p(p^{2}-1)q(q^{2}-1) S_{1}\bigl( \lfloor\tfrac{D}{pq}\rfloor \bigr) + \\
	& + pq(q^{2}-1)\, S_{p} \bigl( \lfloor\tfrac{D}{pq}\rfloor \bigr)+ 
	qp(p^{2}-1)\, S_{q} \bigl( \lfloor\tfrac{D}{pq}\rfloor \bigr)+ 
	pq\, S_{pq} \bigl( \lfloor\tfrac{D}{pq}\rfloor \bigr)
	\,,
	\eas
and the result follows by applying the previous cases.
\end{proof}

We end this section with the two small lemmas used in the last proofs.

\begin{lemma} \label{lem:technicallemma} Let $k$ be a square-free integer. Then
	\[\sum_{m|d} \mu\bigl(\tfrac{d}{m}\bigr) m^{3} e^{*}_{1}(m_{k}^{2}) 
	= \prod_{p|(k,d)}p^{3\nu_{p}(d)-3}(p^{3}-1) \cdot \sum_{m|d_{k}} \mu\bigl(\tfrac{d}{m}\bigr) m^{3} e^{*}_{1}(m^{2}) \]
\end{lemma}

\begin{proof}We first prove the case $k=p$. If $p|d$ then, writing $d=d_{p}\cdot p^{\nu_{p}(d)}$ one gets
	\begin{align*}
	\sum_{m|d} \mu\bigl(\tfrac{d}{m}\bigr) m^{3} e^{*}_{1}(m_{p}^{2}) &= 
	\sum_{m|d_{p}}\sum_{j=0}^{\nu_{p}(d)} \mu\bigl(\tfrac{d}{m p^{j}}\bigr) m^{3}p^{3j} e^{*}_{1}(m^{2}) =\\
	&= p^{3\nu_{p}(d)-3}(p^{3}-1) \sum_{m|d_{p}} \mu\bigl(\tfrac{d}{m}\bigr) m^{3} e^{*}_{1}(m^{2})\,.
	\end{align*}
The result follows from applying the formula recursively on all primes $p$ dividing $k$.
\end{proof}

\begin{lemma} \label{lem:a(d)_and_a(d_p)} Let $p$ be a prime dividing $d$. Then
	\[\a(d)=p^{3\nu_{p}(d)-2} (p^{2}-1)\cdot \a(d_{p})\,.\]
In particular, for $k$ square-free
	\[\a(kd)= \a(d)\cdot\prod_{ \genfrac{}{}{0pt}{}{p|k}{p|d} } p^{3} \cdot \prod_{ \genfrac{}{}{0pt}{}{p|k}{p\nmid d} } p(p^{2}-1)\,.\] 
\end{lemma}

\begin{proof} Write $d=d_{p}\cdot p^{\nu_{p}(d)}$. By the definition of $\a(d)$ and the properties of the M\"{o}bius function one can deduce
	\[\a(d) = d  \!\!\! \sum_{ \genfrac{}{}{0pt}{}{m|d_{p}}{0\le j \le \nu_{p}(d) }} \!\!\! \mu\bigl(\tfrac{d}{m p^{j}}\bigr) m^{2}p^{2j} =
	 d_{p} p^{\nu_{p}(d)} \sum_{m|d_{p}} \mu\bigl(\tfrac{d_{p}}{m}\bigr) m^{2}(p^{2\nu_{p}(d)}- p^{2\nu_{p}(d)-2})\,, \]
and the first formula follows. 

The formula for $\a(kd)$ follows from applying the first formula recursively on primes dividing $k$.
\end{proof}

%% file: sec_volumes.tex
\section{Calculation of the volumes}\label{sec:volumes}

In this section we prove the main theorem of the paper and recalculate the volumes of $\cH(2)$ and of the Prym loci in genus 3 and 4 using our methods.

\subsection{Volume of the gothic locus \texorpdfstring{$\G$}{G}}

In order to calculate the volume of the gothic locus, we write asymptotics for the Euler characteristics of the Teichm\"{u}ller curves $G_{d^{2}}^{r}$ in terms of the value of $d\bmod{6}$.

\begin{lemma}\label{lem:eulergothicasymptotics} The Euler characteristics of the components $G_{d^{2}}^{r}$ are given by
	\[-\chi(G_{d^{2}}^{r}) = 
	\begin{cases}
	\tfrac{13}{720} \cdot \a(d) + \xi(d) & \text{if $(6,d)=1$,} \\
	\tfrac{13}{480} \cdot \a(d) + \xi(d) & \text{if $(6,d)=2$,} \\
	\tfrac{13}{540} \cdot \a(d) + \xi(d) & \text{if $(6,d)=3$,} \\
	\tfrac{13}{360} \cdot \a(d) + \xi(d) & \text{if $(6,d)=6$,}
	\end{cases}
	\]
where the error term $\xi(d)=O(d^{5/2})$ as $d\to\infty$.
\end{lemma}

\begin{proof}First note that $\a(d)=d^{3}\prod_{p|d}(1-1/p^2)$ and therefore by~\eqref{eq:eulercharande} one has
	\[\frac{\chi(X_{d^{2}}(\frakb_{r}))}{d} \le K\cdot \frac{\a(d)}{d} \le K\cdot d^{2}\]
for some constant $K>0$ independent of $d$. In particular, by~\autoref{lem:formulafornonsquare} we have the estimate $\chi(G_{d^{2}}^{r}) \le -\tfrac{3}{2}\chi(X_{d^{2}}(\frakb_{r})) - 2\, \chi(R_{d^{2}}^{r}) + K'\cdot d^{2}$ and the boundary contribution to the Euler characteristic of $G_{d^{2}}^{r}$ is negligible.

By~\autoref{thm:eulergothic} we need to study the asymptotics of $\chi(X_{d^{2}})$ and $\chi(\R[d^{2}])$. Using equations~\eqref{eq:vol16} and~\eqref{eq:eulercharande} and~\autoref{prop:e(d,k)} we have $\chi(X_{d^{2}})=\tfrac{1}{72}\a(d)$ and
	\[2\chi(\R[d^{2}]) = -\frac{1}{3c_{d^{2}}}e(d^{2},6) = 
	\begin{cases}
	-\tfrac{1}{360} \cdot \a(d) + O(d^{5/2}) & \text{if $(6,d)=1$,} \\
	-\tfrac{1}{240} \cdot \a(d) + O(d^{5/2}) & \text{if $(6,d)=2$,} \\
	-\tfrac{1}{270} \cdot \a(d) + O(d^{5/2}) & \text{if $(6,d)=3$,} \\
	-\tfrac{1}{180} \cdot \a(d) + O(d^{5/2}) & \text{if $(6,d)=6$,}
	\end{cases}	\]
and the result follows.
\end{proof}

We finally calculate the volume of the gothic locus.


\begin{proof}[Proof of~\autoref{thm:main}] By~\autoref{prop:volumeformula}, we need to count the number $\cS_{m,m}(\cG)$ of minimal torus covers of fixed degree $m$ in $\G$. 
By~\autoref{thm:areadiscgothic} such number depends on the value of $m$, and is always given by the number of square-tiled surfaces in some gothic Teichm\"{u}ller curves $G_{(m/r)^{2}}^{r}$. By~\autoref{lem:stseulerchar}, this number equals $-6\chi(G_{(m/r)^{2}}^{r})$. 
Altogether, one has
	\bas
	\vol(\G) 
	&=\lim_{D\to\infty} \frac{1}{D^{4}} \sum_{d=1}^{D} \sum_{m|d} \sigma\bigl(\tfrac{d}{m}\bigr) \mid \cS_{m,m}(\cG) \mid = \\
	&= \lim_{D\to\infty} \frac{1}{D^{4}} \sum_{d=1}^{D} \Bigl( \sum_{m|d} \sigma\bigl(\tfrac{d}{m}\bigr)  (-6 \chi(G_{m^{2}}^{1})) + 
	\!\!\! \sum_{\genfrac{}{}{0pt}{}{m|d}{(m,4)=2}} \!\!\!  \sigma\bigl(\tfrac{d}{m}\bigr)  (-6 \chi(G_{(m/2)^{2}}^{2})) + \\
	& 
	\phantom{=} +\!\!\! \sum_{\genfrac{}{}{0pt}{}{m|d}{(m,9)=3}} \!\!\! \sigma\bigl(\tfrac{d}{m}\bigr) (-6 \chi(G_{(m/3)^{2}}^{3})) + 
	\!\!\! \sum_{\genfrac{}{}{0pt}{}{m|d}{(m,36)=6}} \!\!\! \sigma\bigl(\tfrac{d}{m}\bigr) (-6 \chi(G_{(m/6)^{2}}^{6})) \Bigr)\,,
	\eas
where the last equality follows from~\autoref{thm:areadiscgothic}.

The calculations for each of the summands in the limit follow the same lines. We prove in detail the second one and state the rest of the results without proof.

Since the formulae for the Euler characteristics of $G_{d^{2}}^{r}$ depend on the value of $(d,6)$, we first write
	\bas
	\sum_{d=1}^{D} 
	\!\!\! \sum_{\genfrac{}{}{0pt}{}{m|d}{(m,4)=2}} \!\!\!  \sigma\bigl(\tfrac{d}{m}\bigr) &(-6 \chi(G_{(m/2)^{2}}^{2})) = \\
	&= \sum_{d_{0}=1}^{\lfloor D/2 \rfloor} \Bigl(
	\!\!\! \sum_{\genfrac{}{}{0pt}{}{m_{0}|d_{0}}{(m_{0},6)=1}} \!\!\!  \sigma\bigl(\tfrac{d_{0}}{m_{0}}\bigr)  (-6 \chi(G_{m_{0}^{2}}^{2})) \, + \!\!\!\! \sum_{\genfrac{}{}{0pt}{}{m_{0}|d_{0}}{(m_{0},6)=3}} \!\!\!  \sigma\bigl(\tfrac{d_{0}}{m_{0}}\bigr)  (-6 \chi(G_{m_{0}^{2}}^{2})) \Bigr)\,,
	\eas
where we denote $m=2m_{0}$ and $d=2d_{0}$.

Now, using the asymptotics in~\autoref{lem:eulergothicasymptotics} this sum becomes
	\[
	\sum_{d_{0}=1}^{\lfloor D/2 \rfloor } \Bigl(
	\!\!\! \sum_{\genfrac{}{}{0pt}{}{m_{0}|d_{0}}{(m_{0},6)=1}} \!\!\!  \sigma\bigl(\tfrac{d_{0}}{m_{0}}\bigr)  (\tfrac{6\cdot 13}{720}\, \a(m_{0}) + 6\,\xi(m_{0}^{2}) ) + \!\!\!\! \sum_{\genfrac{}{}{0pt}{}{m_{0}|d_{0}}{(m_{0},6)=3}} \!\!\!  \sigma\bigl(\tfrac{d_{0}}{m_{0}}\bigr)  (\tfrac{6\cdot 13}{540}\, \a(m_{0}) + 6\,\xi(m_{0}^{2}) ) \Bigr)\,.
	\]
The error terms can be disregarded using M\"{o}bius inversion, since $\tfrac{1}{D^{4}}\sum d^{5/2} \to 0$. By applying definition~\eqref{eq:sum_convolutions} and an inclusion-exclusion argument one gets
	\begin{multline*}
	\sum_{d=1}^{D} 
	\!\!\! \sum_{\genfrac{}{}{0pt}{}{m|d}{(m,4)=2}} \!\!\!  \sigma\bigl(\tfrac{d}{m}\bigr) (-6 \chi(G_{(m/2)^{2}}^{2})) = \\
	 =\tfrac{13}{360} \Bigl( 3\,S_{1}\bigl(\lfloor \tfrac{D}{2}\rfloor\bigr) - 3\,S_{2}\bigl(\lfloor \tfrac{D}{2}\rfloor\bigr) + S_{3}\bigl(\lfloor \tfrac{D}{2}\rfloor\bigr) - S_{6}\bigl(\lfloor \tfrac{D}{2}\rfloor\bigr) \Bigr) 
	\,.
	\end{multline*}

By~\autoref{prop:S_k} the limit as $D\to\infty$ of this summand divided by $D^{4}$ is
	\[\frac{13\pi^{4}}{360\cdot 2^4}\Bigl(\frac{3}{2^{3}\cdot 3^{2}\cdot 5} - \frac{3}{2^{3}\cdot 3\cdot 5 \cdot 7} + \frac{1}{2\cdot 3^{2}\cdot 5\cdot 13} - \frac{1}{2\cdot 3\cdot 5\cdot 7 \cdot 13} \Bigr) = \frac{43\,\pi^{4}}{2^{8}\cdot 3^{4}\cdot 5^{2} \cdot 7}\,.\]

Proceeding in the same way with the other three summands we get
	\bas
	&\lim_{D\to\infty} \frac{1}{D^{4}} \sum_{d=1}^{D} \sum_{m|d} \sigma\bigl(\tfrac{d}{m}\bigr)  (-6 \chi(G_{m^{2}}^{1})) && \!\!\!\!\!\!= \frac{17\cdot 43\,\pi^{4}}{2^{7}\cdot 3^{4}\cdot 5^{2} \cdot 7} \,, \\
	&\lim_{D\to\infty} \frac{1}{D^{4}} \sum_{d=1}^{D} \!\!\! \sum_{\genfrac{}{}{0pt}{}{m|d}{(m,4)=2}} \!\!\!  \sigma\bigl(\tfrac{d}{m}\bigr)  (-6 \chi(G_{(m/2)^{2}}^{2})) &&\!\!\!\!\!\!= \frac{43\,\pi^{4}}{2^{8}\cdot 3^{4}\cdot 5^{2} \cdot 7}\,, \\ 
	&\lim_{D\to\infty} \frac{1}{D^{4}} \sum_{d=1}^{D} \!\!\! \sum_{\genfrac{}{}{0pt}{}{m|d}{(m,9)=3}} \!\!\! \sigma\bigl(\tfrac{d}{m}\bigr) (-6 \chi(G_{(m/3)^{2}}^{3})) &&\!\!\!\!\!\!= \frac{17\,\pi^{4}}{2^{7}\cdot 3^{5}\cdot 5^{2} \cdot 7} \,, \\
	&\lim_{D\to\infty} \frac{1}{D^{4}} \sum_{d=1}^{D} \!\!\! \sum_{\genfrac{}{}{0pt}{}{m|d}{(m,36)=6}} \!\!\! \sigma\bigl(\tfrac{d}{m}\bigr) (-6 \chi(G_{(m/6)^{2}}^{6})) &&\!\!\!\!\!\!= \frac{\pi^{4}}{2^{8}\cdot 3^{5}\cdot 5^{2} \cdot 7} \,,
	\eas
the sum of which yields the result.

\end{proof}

%

\subsection{Volumes of $\cH(2)$, $\P_{3}$ and $\P_{4}$}

In this section we recalculate the volumes of~$\cH(2)$ and the Prym loci in genus 3 and 4 using our methods.

\begin{theorem}[\cite{ZorVol}]The volume of the minimal stratum in genus~2 is
	\[\vol(\cH(2)) = \frac{\pi^{4}}{2^6 \cdot 3 \cdot 5} 
	\,.\]
\end{theorem}

\begin{proof}Following the same analysis as above and using~\autoref{thm:eulerg2} one has
	\begin{align*}\vol(\cH(2)) 
	&= \lim_{D\to\infty}\frac{1}{D^{4}}\sum_{d=1}^{D}\sum_{m|d} 
	\sigma\bigl(\tfrac{d}{m}\bigr) \bigl(-6 \chi(W_{m^{2}}(2))\bigr)  = \\
	&= \lim_{D\to\infty}\frac{1}{D^{4}}\sum_{d=1}^{D}\sum_{m|d} \sum_{r|m} 
	\tfrac{3}{8}\tfrac{m^{2}(m-2)}{r^{2}} \sigma\bigl(\tfrac{d}{m}\bigr) \mu(r) = \\
	&= \lim_{D\to\infty}\frac{1}{D^{4}}\Bigl( 
	\tfrac{3}{8}S_{1}(D)  - 
	\tfrac{3}{4} \sum_{d=1}^{D}\sum_{\genfrac{}{}{0pt}{}{m|d}{r|m} }
	\tfrac{m^{2}}{r^{2}} \sigma\bigl(\tfrac{d}{m}\bigr) \mu(r) \Bigr)
	\end{align*}

The first summand grows like $\tfrac{\pi^{4}D^{4}}{2^{6}\cdot 3\cdot 5} + O(D^{3})$ by~\autoref{prop:S_k}. One can easily calculate the convolution in the second summand using the same techniques as in~\autoref{sec:numbertheory} to see that it behaves like $\tfrac{D^{3}}{3} \sum_{x=1}^{D}\tfrac{1}{x^{2}} + O(D^{2}\log D)$ for $D\to \infty$, and the result follows.
\end{proof}

The volume of the Prym loci can be calculated from the volumes of the corresponding strata of quadratic differentials. These were calculated by Goujard in~\cite[Appendix A]{goujard} using the convention for the volume form in~\cite{AEZ}:
	\[\vol_{AEZ}(\cQ(-1^{3},3))=\tfrac{5}{9}\pi^{4}\quad\mbox{and}\quad\vol_{AEZ}(\cQ(-1,5))= \tfrac{28}{135}\pi^{4}\,.\]
However, the usual problem with the clash of conventions and the subtleties in the differences between normalisations increase the risk of errors (see~\autoref{subsec:normalisations} for a comparison of normalisations). 

Next we calculate these volumes using our methods. Note that, although strictly speaking the formulae in~\autoref{thm:eulerprym} are only valid for non-square discriminant, one can proceed as in~\autoref{lem:formulafornonsquare} and prove that the boundary contributions are negligible.

\begin{theorem}\label{thm:volPrym3} The volume of the Prym locus in genus~3 is
	\[\vol(\cP_{3}) = \frac{5}{2^8 \cdot 3^3} 
	\, \pi^{4}\,.\]
\end{theorem}

\begin{proof}From~\autoref{thm:areadiscprym} and~\autoref{thm:eulerprym} we get
	\begin{align*}
	\vol(\cP_{3}(4)) 
	&= \lim_{D\to\infty}\frac{1}{D^{4}}\sum_{d=1}^{D} \Bigl(
	\sum_{m|d} -6 \chi\bigl(W^{1}_{m^{2}}(4)\bigr) \sigma\bigl(\tfrac{d}{m}\bigr) + 
	\\
	& \hspace{5.5cm} + \sum_{\genfrac{}{}{0pt}{}{m|d}{(m,4)=2}} \!\!\!\!\! -6 \chi\bigl(W^{2}_{(m/2)^{2}}(4)\bigr) \sigma\bigl(\tfrac{d}{m}\bigr)
	\Bigr) = \\
	&= \lim_{D\to\infty}\frac{1}{D^{4}}\sum_{d=1}^{D} \Bigl(
	15 \sum_{m|d} \sigma\bigl(\tfrac{d}{m}\bigr) \chi(X_{m^{2}}) + 
	\tfrac{15}{2} \!\! \sum_{\genfrac{}{}{0pt}{}{m|d}{(m,2)=2}} \!\! \sigma\bigl(\tfrac{d}{m}\bigr) \chi(X_{m^{2}}) \, + \\
	& \hspace{6cm} + 15 \!\! \sum_{\genfrac{}{}{0pt}{}{m|d}{(m,4)=2}} \!\! \sigma\bigl(\tfrac{d}{m}\bigr) \chi(X_{(m/2)^{2}}) \Bigr) = \\
	& = \lim_{D\to\infty}\frac{1}{D^{4}} \Bigl( \tfrac{5}{24} S_{1}(D) + \tfrac{5}{48} S_{2}(D) + \tfrac{5}{24} \bigl( S_{1}\bigr(\lfloor\tfrac{D}{2}\bigr\rfloor) - S_{2}\bigr(\lfloor \tfrac{D}{2} \rfloor \bigr) \bigr) \Bigr)
	\,,
	\end{align*}
where we have used that
	\[ \sum_{d=1}^{D} \!\! \sum_{\genfrac{}{}{0pt}{}{m|d}{(m,4)=2}} \!\! \sigma\bigl(\tfrac{d}{m}\bigr) \chi(X_{(m/2)^{2}}) = 
	\sum_{d_{0}=1}^{\lfloor D/2 \rfloor} \!\!\! \sum_{\genfrac{}{}{0pt}{}{m_{0}|d_{0}}{(m_{0},2)=1}} \!\!\! \sigma\bigl(\tfrac{d_{0}}{m_{0}}\bigr) \chi(X_{m_{0}^{2}}) =  S_{1}\bigl(\lfloor\tfrac{D}{2}\bigr\rfloor) - S_{2}\bigl(\lfloor \tfrac{D}{2} \rfloor \bigr)\,. \]

Calculating the limits of the four summands by means of~\autoref{prop:S_k} one gets
	\[\vol(\cP_{3}(4)) = \frac{\pi^{4}}{2^{6}\cdot 3^{3}} + \frac{\pi^{4}}{2^{7}\cdot 3^{2}\cdot 7} + \frac{\pi^{4}}{2^{10}\cdot 3^{3}} - \frac{\pi^{4}}{2^{10}\cdot 3^{2}\cdot 7}  =  \frac{5\, \pi^{4}}{2^{8}\cdot 3^{3}}\,. \]

\end{proof}

\begin{theorem}\label{thm:volPrym4} The volume of the Prym locus in genus~4 is
	\[\vol(\cP_{4}) = 
	\frac{7}{2^{9}\cdot 3^{3}\cdot 5} 
	\, \pi^{4}\,.\]
\end{theorem}

\begin{proof}By writing $m=2m_{0}$ and $d=2d_{0}$, the formula for the volume of the Prym locus in genus 4 can be written as 
	\begin{align*}\vol(\cP_{4})
	& = \lim_{D\to\infty}\frac{1}{D^{4}}\sum_{d=1}^{D} \!\!\! \sum_{\genfrac{}{}{0pt}{}{m|d}{m\equiv 0(2)}}\!\!\!  \sigma\bigl(\tfrac{d}{m}\bigr) \bigl(-6 \chi(W_{(m/2)^{2}}(6))\bigr) =  \\
	&= \lim_{D\to\infty}\frac{42}{D^{4}}\sum_{d_{0}=1}^{\lfloor D/2 \rfloor }\sum_{m_{0}|d_{0}} 
	\sigma\bigl(\tfrac{d_{0}}{m_{0}}\bigr) \chi(X_{m_{0}^{2}}) = 
	\lim_{D\to\infty}\frac{7}{12\cdot D^{4}} S_{1}\bigl(\lfloor \tfrac{D}{2} \rfloor\bigr) 
	\,,
	\end{align*}
where the first equality comes from~\autoref{prop:volumeformula} and~\autoref{thm:areadiscprym}, the second one from the formula for the Euler characteristic in~\autoref{thm:eulerprym}, and the final one from~\eqref{eq:eulercharande}. The estimate in~\autoref{prop:S_k} yields the result.
\end{proof}

\subsection{Lattices and volume normalisations for Prym}\label{subsec:normalisations}

The strata $\cQ(\bmb)$ of quadratic differentials are parametrised by the $-1$-eigenspaces $H^{1}_{-}(X,Z(\omega);\bC)$ of the relative cohomology of the canonical double covers $(X,\omega)$. In~\cite{AEZ} and~\cite{goujard} the authors normalise the volume by choosing the following lattice
	\[\Lambda_{AEZ}^{\cQ(\bmb)} \coloneqq \{(Y,\eta)\in H^{1}_{-}(X,Z(\omega);\bC) \,:\, \RPer^{-}(Y,\eta) \subset \bZ\oplus i\bZ \}\,,\]
where $\RPer^{-}(Y,\eta) = \{ \omega(\gamma) \,:\, \gamma \in H_{1}^{-}(X,Z(\omega);\bZ) \}$. 

The Prym loci are a particular case of this construction, that is they are the image under the canonical cover of certain strata of quadratic differentials.

\begin{lemma}\label{lem:prymdoublecover} The canonical double cover of quadratic differentials induces bijections $\cQ(-1^{3},3) \to \cP_{3}\subset \cH(4)$ and $\cQ(-1,5) \to \cP_{4}\subset \cH(6)$.
\end{lemma}

\begin{proof}The maps $(X,\omega)\mapsto (X/J,\omega^{2})$ give the inverse of the canonical double cover map.
\end{proof}

In both cases we are in the minimal stratum, and one can therefore restrict to absolute homology (actually this fact is more general: whenever all zeroes of $\cQ(\bmb)$ have odd order, the antiinvariant space of relative homology $H_{1}^{-}(X,Z(\omega);\bZ)$ agrees with the absolute one $H_{1}^{-}(X;\bZ)$). Moreover, given a surface $(X,\omega)\in\cP_{3}$  (resp. a surface $(X,\omega)\in\cP_{4}$) the polarisation on $H_{1}^{-}(X,Z(\omega);\bZ)$ is of type $(1,2)$ (resp. of type $(2,2)$). This implies the following indices between lattices:
	\be\label{eq:normalisation1_lattice}
	[\Lambda_{abs}^{\cP_{3}}: \Lambda_{AEZ}^{\cQ(-1^{3},3)}] = 2^{4} \qquad\mbox{and}\qquad 
	[\Lambda_{abs}^{\cP_{4}}: \Lambda_{AEZ}^{\cQ(-1,5)}] = 2^{8}\,.
	\ee
To see this, for a $(X,\omega)\in \cP_{3}$ take a basis $H_{1}(X;\bZ)=\langle \alpha_{1},\alpha_{2,1},\alpha_{2,2},\beta_{1},\beta_{2,1},\beta_{2,2}\rangle$ such that $H_{1}^{-}(X;\bZ)=\langle \alpha_{1},\alpha_{2,1}+\alpha_{2,2},\beta_{1},\beta_{2,1}+\beta_{2,2}\rangle$ (cf.~\cite[\S 4]{lanneaunguyen}), and note that $\omega\in\Lambda_{AEZ}^{\cP_{3}}$ are allowed to have half-integral periods $\omega(\alpha_{2,1})=-\omega(\alpha_{2,2})$ and $\omega(\beta_{2,1})=-\omega(\beta_{2,2})$ in $\tfrac{1}{2}(\bZ\oplus i\bZ)$. A similar construction yields the result for $\cP_{4}$ (cf.~\cite[\S 2]{LNPrym4}).

\

There are yet two other sources of disagreement between the AEZ-volumes and ours. The first one, due to another difference in normalisation, is that they define the volume $\vol_{AEZ}$ on $\cQ_{1}(\bmb)$ by disintegration of the Masur-Veech volume element with respect to the area (\cite[\S 4.1]{AEZ}), that is
	\be\label{eq:normalisation2_measure}
	\vol_{AEZ}(B) = 2\dim_{\bC}\cQ(\bmb)\cdot \nu(C(B)) = 2\dim_{\bC}\cQ(\bmb)\cdot \vol(B)\,,
	\ee
where $\vol$ denotes the volume element with our normalisation.

The second one is intrinsic: the poles of the differentials in $\cQ(\bmb)$ are numbered, but they are forgotten under the canonical double cover. This implies (see~\cite[Remark 1.2]{AEZ}) that, for a stratum of quadratic differentials $\cQ(\bmb)$ with $k$ simple poles and its canonical double cover $\cH(\bma)$ one has
	\be\label{eq:normalisation3_poles}
	\vol_{AEZ}(\cQ_{1}(\bmb)) = k! \cdot \vol_{AEZ}(\cH_{1}(\bma))\,.
	\ee

All factors considered, one has (cf.~\cite[Appendix A]{goujard}):
	\[\vol_{AEZ}(\cQ_{1}(-1^{3},3)) = \frac{5}{9}\pi^{4} = 2^{4}\cdot 2^{3}\cdot 3! \cdot \vol(\cP_{3})\,,\]
	where the $2^{4}$ corresponds to the different lattice normalisation in~\eqref{eq:normalisation1_lattice}, the $2^{3}$ to the measure normalisation in~\eqref{eq:normalisation2_measure} and the $3!$ to the numbering of the poles in~\eqref{eq:normalisation3_poles}, and
	\[\vol_{AEZ}(\cQ_{1}(-1,5)) = \frac{28}{135}\pi^{4} = 2^{8}\cdot 2^{3}\cdot \vol(\cP_{4})\,,\]
	where the $2^{8}$ corresponds to the different lattice normalisation in~\eqref{eq:normalisation1_lattice} and the $2^{3}$ to the measure normalisation in~\eqref{eq:normalisation2_measure}.

%% file: sources.bbl
\newcommand{\etalchar}[1]{$^{#1}$}